\documentclass[11pt]{article}
\newif\iffigures\figurestrue % Set to false to produce document without
\newif\iftables\tablestrue % Set to false to produce document without
\iffigures\usepackage{curves}\fi
\usepackage{amssymb}
\newcommand{\df}{\displaystyle\frac}
\newcommand{\qed}{\rule{3mm}{3mm}}
\newcommand{\dd}{{, \ldots ,}}

\newcommand{\hsp}{\hspace*{\parindent}}

\newcommand{\Ga}{\Gamma}
\newcommand{\ga}{\gamma}

\newcommand{\In}{\infty}

\newcommand{\btheta}{{\mbox{\boldmath $\theta$}}}

\newcommand{\bga}{{\mbox{\boldmath $\gamma$}}}
\newcommand{\RR}{{\mathbb R}}
\newcommand{\ZZ}{{\mathbb Z}}

\newcommand{\QQ}{{\mathbb Q}}
\newcommand{\CC}{{\mathbb C}}

\newcommand{\sF}{{\cal F}}

\newcommand{\sM}{{\cal M}}
\newcommand{\sP}{{\cal P}}

\newcommand{\sS}{{\cal S}}

\newcommand{\bx}{{\bf x}}

\newcommand{\bz}{{\bf z}}
\newcommand{\bv}{{\bf v}}

\font\phvr=phvr at 11pt
\newcommand{\hc}[1]{\mbox{\phvr  #1}}
\newcommand{\beql}[1]{\begin{equation}\label{#1}}
\newcommand{\eeq}{\end{equation}}
\newcommand{\BE}{\begin{equation}}
\newcommand{\EE}{\end{equation}}
\newcommand{\numb}[1]{(\ref{#1})} % refer to equation number
\newtheorem{thm}{Theorem}[section]
\newtheorem{cor}{Corollary}[section]
\newtheorem{lemma}{Lemma}[section]
\newtheorem{definition}{Definition}[section]
\newenvironment{defn}{\begin{definition}\em}{\end{definition}}
\newcommand{\twobytwo}[4] % 2 by 2 matrix
   {\left[ \begin{array}{cc} #1 & #2 \\ #3 & #4 \end{array} \right]}
\newcommand{\hpl}{{\mathfrak H}} % Notation used for hyperbolic plane
\newcommand{\thfl}{\floor {\theta + \frac 1 2 }}

\newcommand{\floor}[1]{\left\lfloor #1 \right\rfloor}
\newcommand{\app}[1]{p_{#1}-q_{#1}\theta} % Error of approximation
\newcommand{\papp}[1]{\left(\app{#1}\right)} % Error in parentheses
\newcommand{\nrapp}[1] % Error with primes and superscript n
   {p_{#1}^{\prime(n)}-q_{#1}^{\prime(n)}\theta}
\newcommand{\interval}{\mathop{\rm int}}
\font\phvrsm=phvr at 8pt %\phvr small, for figures
\font\phvrvs=phvr at 5.5pt %\phvr at half size
\newcommand{\hcsm}[1]{\mbox{\phvrsm #1}}
\newcommand{\hcvs}[1]{\mbox{\phvrvs #1}}
\newcommand{\clap}[1]{\hbox to 0pt{\hss #1\hss}}
\newcommand{\vclap}[1]{\vbox to 0pt{\vss #1\vss}}
\newcommand{\smallmat}[4]{\left[{#1 \atop #3}{#2 \atop #4}\right]}
\newcommand{\un}{\underline 1}
\newcommand{\dblinf}[2] % Doubly-infinite sequence, with the break
   {$[\,\ldots,#1,$ & $\,#2,\ldots\,]$}
\newcommand{\dblunf}[2] % Doubly-infinite sequence, with the underlined
   {\dblinf{#1}{\un,#2}}
\newcommand{\symrow}[5] % Symmetric row, containing a 5-term sequence
   {\dblinf{#1,#2}{#3,#4,#5}&\dblinf{#5,#4}{#3,#2,#1}\\}
\newcommand{\dblfin}[2] % Doubly-finite sequence, with the break between
   {$[#1,$ & $\,#2]$}
\newcommand{\dblcen}[2] % Doubly-finite sequence, with the 1_c lined up
   {\dblfin{#1}{1_c,#2}}
\newcommand{\finrow}[8] % Four doubly-finite sequences in one row
   {\dblcen{#1}{#2}&\dblcen{#3}{#4}&\dblcen{#5}{#6}&\dblcen{#7}{#8}\\}
\newcommand{\shortrow}[4] % \shortrow split in two
{\dblcen{#1}{#3}&\dblcen{#1}{#4}\\
 \dblcen{#2}{#3}&\dblcen{#2}{#4}\\ \hline}
\newcommand{\sglinf}[2] % Singly-infinite sequence, with the break
   {$[a_0;#1,$ & $\,#2,\ldots\,]$}

\catcode`\@=11
\renewcommand{\section}{
        \setcounter{equation}{0}
	\setcounter{figure}{0}
	\setcounter{table}{0}
        \@startsection {section}{1}{\z@}{-3.5ex plus -1ex minus
        -.2ex}{2.3ex plus .2ex}{\large\bf}
        }
\catcode`\@=12
\makeatletter
\def\eqalignno#1{\displ@y \tabskip\@centering
  \halign to\displaywidth{\hfil$\@lign\displaystyle{##}$\tabskip\z@skip
    &$\@lign\displaystyle{{}##}$\hfil\tabskip\@centering
    &\llap{$\@lign##$}\tabskip\z@skip\crcr
    #1\crcr}}
\makeatother

\makeatletter
\def\@sect#1#2#3#4#5#6[#7]#8{\ifnum #2>\c@secnumdepth
     \def\@svsec{}\else
     \refstepcounter{#1}\edef\@svsec{\csname the#1\endcsname.\hskip .75em }\fi
     \@tempskipa #5\relax
      \ifdim \@tempskipa>\z@
        \begingroup #6\relax
          \@hangfrom{\hskip #3\relax\@svsec}{\interlinepenalty \@M #8\par}%
        \endgroup
       \csname #1mark\endcsname{#7}\addcontentsline
         {toc}{#1}{\ifnum #2>\c@secnumdepth \else
                      \protect\numberline{\csname the#1\endcsname}\fi
                    #7}\else
        \def\@svsechd{#6\hskip #3\@svsec #8\csname #1mark\endcsname
                      {#7}\addcontentsline
                           {toc}{#1}{\ifnum #2>\c@secnumdepth \else
                             \protect\numberline{\csname the#1\endcsname}\fi
                       #7}}\fi
     \@xsect{#5}}
\def\@begintheorem#1#2{\it \trivlist \item[\hskip \labelsep{\bf #1\ #2.}]}
\makeatother
\begin{document}
\thispagestyle{empty}
\begin{center}
{\Large {\bf Cutting Sequences for Geodesic Flow on the \smallskip \\
Modular Surface and Continued Fractions }} \\
\vspace{\baselineskip}
{\large {\em David J. Grabiner}}\footnote{Supported by an NSF Postdoctoral
Fellowship.  Some of this work was done while visiting the Mathematical
Sciences Research Institute; MSRI is supported by NSF Grant
DMS-9022140.} \\
  \vspace{.5\baselineskip}
Department of Mathematics\\
Arizona State University\\
Tempe, AZ 85287-1804\\
\vspace{.5\baselineskip}
{\large {\em Jeffrey C. Lagarias}}\footnote{Some of this work was done
while visiting the Mathematical Sciences Research Institute.
} \\
\vspace{.5\baselineskip}
AT\&T Labs \\
Florham Park, NJ 07932-0971 \\
\vspace{1\baselineskip}
{April 1, 2001}\\
\vspace{2\baselineskip}
{\em Abstract}
\end{center}
This paper describes the cutting sequences of geodesic flow on the
modular surface $\hpl / PSL ( 2, \ZZ )$ with respect to the standard
fundamental domain $\sF = \{ z = x + iy : $ $ - \frac{1}{2} \leq x \leq
\frac{1}{2}$ and $| z | \geq 1 \}$ of $PSL(2, \ZZ )$.  The cutting
sequence for a vertical geodesic $\{ \theta + it : t > 0 \}$ is
related to a one-dimensional continued fraction expansion for $\theta$,
called the one-dimensional Minkowski geodesic continued fraction (MGCF)
expansion, which is associated to a parametrized family of reduced bases
of a family of 2-dimensional lattices.  The set of cutting sequences 
for all geodesics
forms a two-sided shift in a symbol space $\{ \bar{\hc{L}} ,
\bar{\hc{R}} , \bar{\hc{J}} \}$ which has the same set of forbidden
blocks as for vertical geodesics.  We show that this shift is not a
sofic shift, and that it characterizes the fundamental domain $\sF$ up
to an isometry of the hyperbolic plane $\hpl$.
We give conversion methods between the cutting sequence for the
vertical geodesic $\{\theta + it: t > 0\}$, the MGCF expansion of $\theta$
and the additive ordinary continued fraction (ACF) expansion of $\theta.$
We show that the cutting sequence and MGCF expansions can each
be computed from the other by a finite automaton, and the ACF
expansion of $\theta$ can be computed from
the cutting sequence for the vertical geodesic $\theta + it$
by a finite automaton.  However, 
the cutting sequence for a vertical geodesic
cannot be computed from the ACF expansion
by any finite automaton, but there is an algorithm to compute
its first $\ell$ symbols when given as input  the first $O(\ell)$ symbols
of the ACF expansion, which takes time $O(\ell^2)$ and
space $O(\ell)$.

\setcounter{page}{0}
\setcounter{footnote}{0}
\newpage
\begin{center}
{\Large {\bf Cutting Sequences for Geodesic Flow on the \smallskip \\
Modular Surface and Continued Fractions }} \\
\vspace{\baselineskip}
{\large {\em David J. Grabiner}} \\
\vspace{.5\baselineskip}
Department of Mathematics\\
Arizona State University\\
Tempe, AZ 85287-1804\\
\vspace{.5\baselineskip}
{\large {\em Jeffrey C. Lagarias}} \\
\vspace{.5\baselineskip}
AT\&T Labs \\
Florham Park, NJ 07932-0971 \\
\vspace{1\baselineskip}
% (\today)
\vspace{2\baselineskip}
\end{center}
\section{Introduction}
\hsp
This paper describes the symbolic dynamics of cutting sequences of
geodesics for the standard
fundamental domain $\sF =  \{ z =  x  +  iy  :  - \frac{1}{2}  \leq  x <
\frac{1}{2}$,
$|z|  \geq 1  \}$ of the modular group
$PSL(2, \ZZ )$ acting on the upper half plane $\hpl =  \{ z \in \CC  :  Im
(z)  >  0  \}$.

We study in particular the cutting sequences of vertical geodesics
$\{  \theta  +  it  :  t > 0  \}$ for $\theta  \in  \RR$.
These particular geodesics
are related to a continued fraction expansion
introduced in 1850 by Hermite \cite{Her51} in terms of quadratic forms, and
studied by
Humbert \cite{Hum16a,Hum16b}.
Our motivation for studying them was that
they appear in the 
 one-dimensional case of a 
multidimensional continued fraction 
introduced in  Lagarias \cite{Lag94}. This
expansion, called the {\em Minkowski geodesic continued fraction expansion}
 (MGCF expansion), is based on following a parametrized family of 
lattice bases in 
$GL( d +1, \ZZ ) \backslash GL( d + 1, \RR )$ as the parameter $t$ varies; 
in the one-dimensional case the family associated to $\theta$ is 
$\left[
\begin{array}{ll} 1 & 0
\\ - \theta & t
\end{array}
\right]$.
The MGCF expansion has the merit that it finds 
 good Diophantine approximations
in all dimensions. However, the
one-dimensional case of this expansion
does not coincide with the ordinary continued
fraction expansion.  We show in \S 3.4 that the one-dimensional MGCF
expansion of $\theta$ is
essentially equivalent to the cutting sequence expansion of
the vertical geodesic $\{ \theta + it : 0 < t < \infty \}$.
We use this connection to determine the precise
 relation of the MGCF expansion 
of $\theta$ to the additive
continued fraction expansion of $\theta$; the
{\em additive continued fraction} is the variant of the
ordinary continued fraction which includes all intermediate
convergents. 

A second  motivation for studying cutting sequences on 
$\hpl / PSL (2, \ZZ)$ arises from their use as symbolic codings of
geodesics in the study of geodesic flow on constant negative
curvature surfaces of finite volume; we give more background on this
at the end of the introduction. 
Cutting sequence encodings are of special interest 
because they  apparently
encode more information about the geodesic flow on a Riemann surface
than other symbolic encodings of geodesics, which 
have a simple description (shift of finite type) but
only retain topological and not conformal information
about the Riemann surface.
Adler and Flatto \cite[\S 10]{AdlFla91} raised the question
whether cutting sequence encodings preserve conformal information about
the Riemann surface and also retain information
about the lines on the Riemann surface
used to define the cuts, i.e. the shape of the fundamental domain in
the universal cover.  Consider any
hyperbolic polygon $\sP$ in $\hpl$
which is a fundamental domain for a properly discontinuous group
$\Gamma$ acting on $\hpl$, such that $\hpl / \Gamma$ has finite
volume, and suppose that $\sP$ has $| \sP | $ sides.  Let $\Sigma_\sP$
denote the closed subshift of the shift on $| \sP | $ letters generated
by the cutting sequences for $\sP$ for geodesic flow on general $\hpl /
\Gamma$, as defined in \S 2.2.
Adler and Flatto \cite[p. 300] {AdlFla91}\footnote{They raise the question for the $(8g-4)$-sided
fundamental polygons discussed in \cite{AdlFla91}, for genus $g \geq 2$.} ask:
does $\Sigma_\sP$
determine $\sP$ up to conformal isometry and $\Gamma$ up to hyperbolic
conjugacy?  We show here that this is the case for the standard
fundamental domain $\sF$ of $PSL(2, \ZZ)$.

Our main results on cutting sequence encodings are as follows.
\begin{description}
\item{(1)}
The set of cutting sequence expansions for (irrational) vertical geodesics
$\Pi_\sF^0$ is characterized in terms of
forbidden blocks (Theorem~6.1).
A generating set of forbidden blocks is enumerated.
(Theorems~\ref{1h1mthm} and \ref{blockthm}.).
The number of minimal forbidden blocks of length at most $k$ grows
exponentially in $k$.
(Theorem~6.3).
\item{(2)} 
The cutting sequence expansions for all geodesics
comprise
(essentially) a closed two-sided subshift $\Sigma_\sF$ of the full shift on
three symbols $\{ \bar{\hc{L}}, \bar{\hc{R}}, \bar{\hc{J}} \}$.  The set
of forbidden blocks of $\Sigma_\sF$ is the same as that for $\Pi_\sF^0$.
(Theorem \ref{twosidethm}).
Each symbol sequence in $\sum_\sF$ corresponds to a unique
oriented geodesic on $\hpl / PSL (2, \ZZ)$, with the exception of the two
sequences $\bar{\hc{L}}^\In$ and $\bar{\hc{R}}^\In$.
In the converse direction, each oriented geodesic gives rise to a finite
number of
shift-equivalence classes of
symbol sequences in $\Sigma_\sF$.
This number is at most eight if it is not a periodic geodesic.
(Theorem~7.2).
The shift $\Sigma_\sF$ is
not a sofic system. (Theorem \ref{notsofic})
\item{(3)} 
The shift $\Sigma_\sF$ characterizes $\sF$ up to an isometry
of $\hpl$.
If $\sP$ is a hyperbolic polygon of finite area (possibly
with some ideal vertices) which is a fundamental domain of a discrete
subgroup $\Gamma$ of $PSL (2, \RR )$ and if the subshift $\Sigma_\sP$ is
isomorphic to $\Sigma_\sF$ by a permutation of symbols, then there is a
hyperbolic isometry $g \in PSL(2, \RR)$ such that $\sF = g \sP$ and $PSL
(2, \ZZ) = g \Gamma g^{-1}$.  (Theorem \ref{uniquepoly}).
\end{description}

In obtaining result (1) we show in \S 5
that vertical geodesics have special features (not shared by general
geodesics) that facilitate characterizing forbidden blocks: any
vertical geodesic that hits a $PSL(2,\ZZ)$-translate of a corner of the
fundamental domain must have rational $\theta$,
and
each such geodesic hits at most one such corner,
with the exception of those $\theta \equiv \frac{1}{2} ( \bmod ~1)$,
which hits two corners.
Furthermore,
the continued fraction expansion of any rational $\theta$ such that $\{
\theta + it : t > 0\}$ hits a corner
has a symmetry which can be described in terms
of the linear fractional transformation $\bz \rightarrow \hc{N} \bz$ with
$\hc{N} = \left[ \begin{array}{ll}
1 & 2 \\
2 & 1
\end{array}
\right]$; see Theorem~5.1.
In formulating result (3)
we {\em define} cutting sequence expansions for
geodesics that hit a corner of the fundamental domain $\sF$ of $\hpl /
PSL (2, \ZZ )$ to be limits of general-position geodesics that hit no
corner.
This procedure assigns infinitely many different cutting sequences
to certain periodic geodesics.
We interpret result (3) as showing that the shift $\Sigma_\sF$
encodes conformal information about $\sF$.

As a preliminary to obtaining the above results, in \S2--\S5 we determine
precise relations of the cutting sequence
of the geodesic $\theta + it$ 
to the Minkowski geodesic continued fraction expansion of $\theta$
and to the additive version of the
 ordinary continued fraction expansion of $\theta,$ 
which we call the ACF expansion. (The  ACF expansion, which  is described
in \S 3.1, includes 
all intermediate convergents.)
These give algorithms for converting between the symbolic expansions,
and obtain results on how hard it is computationally to convert 
between these different
symbolic expansions, as follows.

\begin{description}
\item{(1)} 
The cutting sequence expansion of $\{\theta + it: t > 0\}$ and the
Minkowski geometric continued fraction expansion of $\theta$ can
each be computed from the other by a finite 
automaton. (Theorem~\ref{MGCFtoCut}).
The additive ordinary continued fraction expansion of
$\theta$ can be computed from either of these expansions
by a finite automaton. (Theorems \ref{MGCFtoOCF} and \ref{CStoOCF}).

\item{(2)}
The Minkowski geodesic continued fraction expansion of $\theta$
{\em cannot} be computed from the additive ordinary continued
fraction expansion of $\theta$ by a finite automaton. 
(Theorem \ref{noautthm}).

\item{(3)} 
The  Minkowski geodesic continued fraction expansion
of $\theta$ can be computed from the additive ordinary continued
fraction expansion of $\theta$ in quadratic time and linear space
on a random access machine,
i.e. the first  $l$ $MGCF$ digits can be computed in time $O(\ell^2)$ and
space $O(\ell)$, from  $O(\ell)$ symbols of the ACF expansion. 
(Theorem~\ref{turing}).
\end{description}

The notion of finite automaton (finite state machine) is 
described in \S 3.5.
The key feature of a finite automaton is that it has
a fixed finite amount of memory. Result (1) requires that
given the successive symbols in the input  expansion  
there is a an absolute bound $B$ so that after scanning $B$ 
consecutive input symbols at least one new output symbol is
determined. This property may be called bounded look-ahead.
Result (2) comes from
the fact that the amount of look-ahead needed to compute one
MGCF symbol from the ordinary continued fraction 
expansion can be unbounded. Result (3) is based on
a conversion algorithm given in Theorem~\ref{1h1mthm}. To compute
a given output symbol it potentially requires remembering
the entire input string of symbols to that point,
which is the source of the linear space requirement. 
This computation can be
carried out in the given time and space in the 
random access
machine(RAM) model for computation, which
allows unbounded storage.
It can also be carried out on a one-tape Turing machine 
using  $O(\ell)$ space
with a time bound polynomial in $\ell$ on a one-tape Turing
machine.  For the random access machine (RAM) and
Turing machine models
of computation, see Aho, Hopcroft and Ullman ~\cite{AhoHopUll74}.
Results (1)-(3)
delineate the computational relations between
the ordinary continued fraction expansion, the MGCF expansion
and the cutting sequence expansion. Result (2), which
is the significant result here, shows that 
cutting sequences encode information more concisely 
than the additive ordinary continued fraction expansion.

To conclude this introduction, and to put the results
above in a more general context, we recall background on symbolic
codings of geodesic flow on a compact Riemann surface with 
a finite number of punctures.
If $\Gamma$ is a finitely
generated discrete group of conformal isometries of the hyperbolic plane
$\hpl$, such that $\hpl / \Gamma$ has finite volume, then $\hpl / \Gamma$
is a compact Riemann surface\footnote{There are some mild extra
conditions needed on $\Gamma$ for compactness, see Lehner
\cite[pp.~203--205]{Leh64} for sufficient conditions.  It suffices for
$\Gamma$ to have a fundamental region $\sF$ that is a hyperbolic polygon
with a finite number of sides, with only parabolic vertices at
infinity.}  minus a finite number of punctures, and geodesic flow on
$\hpl / \Gamma$ contains both topological and conformal information about
this Riemann surface.  Symbolic codings of geodesics were introduced by
Hadamard \cite{Had98} in 1898 as a way of understanding the complicated
motions on geodesics on such surfaces.  E. Artin \cite{Art24} showed in
1924 that there exist dense geodesics on the modular surface $\hpl /PSL
(2, \ZZ )$, using symbolic encodings with continued fractions, and in
1935 G. Hedlund \cite{Hed35} used Artin's coding to show that geodesic
flow on this surface was ergodic.  Some other symbolic encodings
were ``boundary expansions'' by Nielsen \cite{Nei27} and
cutting sequences with respect to a fundamental domain of $\Gamma$ by
Koebe \cite{Koe29} and Morse~\cite{Mor66}.  Cutting sequence expansions
can be viewed as a generalization of the reduction theory of indefinite
binary quadratic forms to arbitrary finitely-generated Fuchsian groups,
see Katok \cite{Kat86}.  More recently Bowen and Series \cite{BowSer79}
and Series \cite{Ser85a,Ser85b} used ``modified boundary
expansions'' to obtain a particularly simple symbolic expansion from
which strong forms of ergodicity could be deduced; their codings
are sofic systems.
Adler and Flatto
\cite{AdlFla82,AdlFla91} used ``rectilinear map'' codings of
geodesic flow on $\hpl / \Gamma$ to explain why certain specific maps of
the interval, such as the continued fraction map and backwards continued
fraction map, have invariant measures of a simple form.  They showed
that these invariant measures are inherited from the invariant measure
on geodesic flow (Liouville measure) using a cross-section map followed
by a factor map.  The codings of 
Adler and Flatto \cite{AdlFla91}
 have a particularly 
simple structure -- they are shifts of finite type -- and they preserve
topological information about the Riemann surface, but they lose
conformal information, for they are identical for all Riemann surfaces
of the same genus $g \ge 2$.  (See \cite[Theorem 8.4 and
Section 10]{AdlFla91}.) 
In this respect
it is of interest whether cutting sequence encodings preserve conformal
information. Besides all these encodings, two other 
encodings of geodesic flow on the modular surface 
$\hpl/ PSL (2, \ZZ )$  with interesting Diophantine approximation properties
appear in Arnoux \cite{Arn94} and Lagarias and Pollington \cite{LagPol95}.

With regard to this general framework, 
we note that the proofs in this paper
heavily rely on specific facts about $ PSL (2, \ZZ )$. 
The methods used  may  extend in some form
to congruence  subgroups of $ PSL (2, \ZZ )$, i.e. arithmetic
Fuchsian groups, but
do not appear likely to apply to general  finitely generated Fuchsian
groups $\Gamma$. 

For general terminology in symbolic dynamics we follow Adler and
Flatto~\cite{AdlFla91} and Lind and Marcus~\cite{Lin95}.
{}From the viewpoint of symbolic dynamics, the
associated shift $\Sigma_\sF$ studied in this paper is
an example of a naturally occurring shift more
complicated than the ones currently having a
well-developed theory, cf.\cite{Lin95}.

\section{Cutting Sequences}
\hsp
We describe cutting sequences for finite-sided hyperbolic polygons
that are the fundamental domain of a given discrete subgroup $\Gamma$ of
$PSL (2, \RR )$ that acts discontinuously on $\hpl$,
and then specialize to the standard fundamental domain $\sF$ of the
modular group $PSL (2, \ZZ )$.
General references for geodesic flow and for cutting sequences are
Adler and Flatto \cite{AdlFla91} and Series~\cite{Ser91}.

\subsection{Cutting Sequence Shifts}
\hsp Let $\hpl$ denote the hyperbolic plane, represented as the upper
half-plane $\hpl = \{z = x+ iy : y > 0 \}$ with hyperbolic line element
$ds^2 = \frac{1}{y^2} (dx^2 + dy^2 )$ and volume $\frac{dxdy}{y^2}$.  An
oriented geodesic $\bga = \langle \theta_1 , \theta_2 \rangle$ in
$\hpl$ is uniquely specified by its two ideal endpoints $\theta_1 ,
\theta_2 \in \RR \cup \{ \In \}$, with $\theta_1 \neq \theta_2$.  It is
oriented from its {\em head} $\theta_1$ to its {\em foot} $\theta_2$.  A
{\em vertical geodesic} $\langle \In , \theta \rangle$ is the vertical
line $\theta + it$, with $t$ decreasing from $\In$ to 0.  Other
geodesics are semicircles perpendicular to the real axis with endpoints
$\theta_1$ and $\theta_2$, centered at $\frac{\theta_1 + \theta_2}{2}$,
and oriented to start at $\theta_1$ and travel to $\theta_2$.

Let $\Gamma$ be a finitely generated discrete subgroup of $PSL(2, \RR )$
which acts properly discontinuously on the upper half-plane $\hpl = \{
x+iy : y > 0 \}$ as a group of linear fractional transformations.  We
suppose that $\hpl / \Gamma$ has finite volume with respect to the
hyperbolic volume on $\hpl$.  Let $\sP$ be any finite-sided convex
hyperbolic polygon which is a {\em fundamental domain} for $\Gamma$,
i.e. $\sP$ is such that the set of $\Gamma$-translates $\{g \sP : g \in
\Gamma \}$ tessellates $\hpl$, and no element $g \in \Gamma$ sends
$\sP$ to itself except the identity.  We allow polygons $\sP$ to have
one or more ideal vertices (``cusps''), and they may also have extra
vertices, called {\em elliptic vertices}, which have an angle of $\pi$
and which fall in the middle of a side of $\sP$ regarded as a geometric
object, as discussed below.  Fundamental domains exist for each such
finitely generated $\Gamma$, cf. [4, Theorem~10.1.4].  Such polygons
$\sP$ necessarily have an even number of ``sides,'' where by convention
each elliptic vertex divides a side of the geometric object $\sP$ into
two ``sides.''  If $\sP$ is to be a fundamental domain for $\Gamma$, the
``sides'' of $\sP$ must be identified under the action of $\Gamma$.
Specifically, there is a pairing of the ``sides'' of $\sP$ which assigns
to each ``side'' $s$ of $\sP$ a ``side'' $s'$ and an isometry $g(s,s')
\in \Gamma$ such that $g(s,s')$ maps $s$ to $s'$.  Furthermore $(s')' =
s$ and
\[
g(s',s) = (g(s,s'))^{-1},
\]
while if $s = s'$ then $g(s,s)$ is
the identity on $s$.  The pairing indicates how $\sP$ tessellates
$\hpl$.  For $g=g(s,s')$, the domain of $\sP$ is adjacent to $\sP$
in such a way that the ``side'' labeled $s$ of $g\sP$ coincides with
the ``side'' labeled $s'$ of $\sP$, while $g(Int(\sP))\cap
Int(\sP)=\emptyset$.
If $s = s'$ then
necessarily $g^2$ is the identity on $\hpl$, and $g$ is a hyperbolic
reflection.  The {\em neighborhood set} of $\sP$ \beql{eq21} N_\Ga ( \sP
) : = \{g \in \Gamma : g = g (s,s') ~~ \mbox{for ~some~side~$s$ of $\sP
\}$} ~.  \eeq is a set of generators of the group $\Gamma$.  It will
comprise the set of symbols used in the cutting sequence encodings of
geodesics described below.

The number of elements of $N_\Gamma ( \sP )$ is generally equal to the
number of ``sides'' of $\sP$, but can be less if there are two or
more elements $g(s, s')$ that are equal.
For convex polygons $\sP$ this situation can only arise when there are
hyperbolic reflections $g = g(s,s')$ with $s \neq s'$,
in which case $g(s,s') = g(s',s)$.
In such a case the ``sides'' $s$ and $s'$ must lie on a common geodesic
of $\hpl$ and share a common endpoint,\footnote{For non-convex $\sP$ this
may fail to hold.}
which is then a vertex of $\sP$ having an angle of $\pi$.
We call any such vertex an {\em elliptic vertex}.
These two ``sides'' $s$ and $s'$ are then assigned the same
generator, so if we erase the elliptic vertex and glue them together, we
get a single geometric side of the polygon $\sP$ regarded as a
geometric object.
We therefore have:
For a convex fundamental domain $\sP$ of a finitely generated group
$\Gamma$, the number of
elements of $N_\Gamma (P)$ is equal to the number of geometric sides of
$\sP$.  In particular, {\em the number of elements of $N_\Gamma (P)$ may be
odd}.

As an example, the standard fundamental domain $\sF $ of $\Gamma = PSL (2,
\ZZ )$,
is geometrically
a triangle with one ideal vertex, but as a fundamental domain
for $\Gamma$ is a quadrilateral having an elliptic vertex at $z = i$.
In this case $N_\Gamma ( \sF )$ has three elements.

We associate to geodesic flow on $\hpl/ \Gamma$ the (two-sided) {\em
cutting sequence subshift} $\Sigma_{\sP , \Gamma}$ which is a closed
subshift of the full shift on $|N_\Gamma ( \sP )|$ letters.

\begin{defn}
The set $G_{\sP,\Gamma}^0$ of {\em general position geodesics} for
$\sP$ consists of those geodesics
$\gamma$ on $\hpl$
which intersect the interior of $\sP$, and which are transverse to all
$\{g ( \sP ) : g \in \Gamma \}$ in the sense that $\gamma$ contain no
vertex nor part of any side of positive length of any translated polygon
$\{g \sP : g \in \Gamma \}$.  If $\sP$ contains ideal
vertices (``cusps''), then the endpoints of $\gamma$ must not coincide with any
$\Gamma$-translate of any such vertex.
\end{defn}

To each $\gamma \in G_{\sP,\Gamma}^0$,
we assign the doubly-infinite {\em cutting sequence}
\beql{eq22}
C( \ga ) : = ( \ldots , g_{-2} , g_{-1}, g_0 , g_1 , g_2 , \ldots )
\eeq
in which all $g_i \in N_\Gamma( \sP )$, as follows.  We first label each
``side'' $s$ of the fundamental domain $\sP$ with the label $g(s,s')$ on
its ``inside'' edge in $\sP$.  We similarly label all the sides of the
translated domain $\{h \sP : h \in \Ga \}$, so that the edge of $h\sP$
which is $hs$ has the label $g(s,s')$.  Every ``side'' is now
assigned two labels, on its ``inside'' edge and ``outside'' edge,
because it abuts two copies of a fundamental domain.  (Figure~2.1 below
indicates such labels.)  The geodesic $\ga = \langle \theta_1 , \theta_2
\rangle$ visits in order a sequence of translates $\{h_i \sP ;
i \in \ZZ \}$ in which $h_i \in \Gamma$ and $h_0$ is the identity.
At a crossing of the edge from $h_{i-1}  \sP$ to $h_i \sP$ we assign the
symbol $g_{i}$ on the side of the domain that
the geodesic {\em enters}.
(This is the convention of Katok \cite{Kat96};
and opposite to that of Adler and Flatto [3, p.~243], who
use the symbol $g_{i}$ on the exit edge.)
This convention yields:

\begin{lemma}\label{translemma}
The cutting sequence $(g_1,\ldots,g_j)$ follows a geodesic from $\sP$
to the translated domain $\sP_j=h_j\sP$ with
\BE
h_j=g_1 g_2 \cdots g_j.
\EE
\end{lemma}

{\bf Proof.}
The edge $h_i s$ of $h_i \sP$ that borders $h_{i-1} \sP$ is assigned the edge
label inside $h_i \sP$ of $g_i = g(s,s')$.
Now we have
$(h_{i-1} g(s,s') h_{i-1}^{-1})h_{i-1}\sP=h_i\sP$, hence
\beql{eq23}
h_{i-1} g_i = h_i ~,
\eeq
see Figure~2.1.
$\qed$
\iffigures
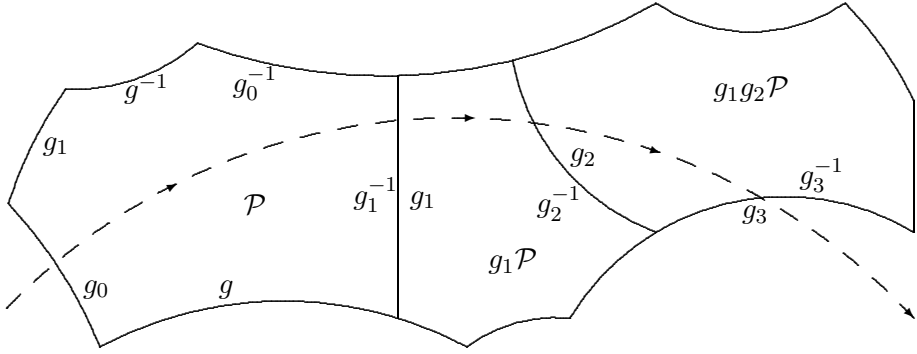
\begin{figure}
%Figure 2.1
\setlength{\unitlength}{.3in}
% This figure is based on a sketch rather than any actual values; angles
% and arcs are chosen so that paths close up properly
\begin{picture}(18,7)(0,1)
% Arcs bounding outer domains
% endpoints 9.2,1.5 and 2.8,1.5, passes through 8,2
\put(6.000,-4.490){\arc(3.200,5.990){56.224}}
% endpoints 2.8,1.5 and 1.2,4
\put(-4.688,-1.530){\arc(7.488,3.030){21.175}}
% endpoints 2.2,6 and 1.2,4
\put(7.900,1.900){\arc(-5.700,4.100){18.325}}
% endpoints 2.2,6 and 4.5,6.8
\put(2.370,9.216){\arc(-0.170,-3.215){44.435}}
% endpoints 4.5,6.8 and 12.5,7.5, passes through 8,6.230 and 10,6.5
\put(7.747,15.758){\arc(-3.247,-8.958){49.847}}
% endpoints 12.5,7.5 and 15.8,7.5
\put(14.150,9.950){\arc(-1.650,-2.450){67.918}}
% endpoints 17,5.8 and 15.8,7.5
\put(11.394,3.117){\arc(5.606,2.683){19.275}}
% endpoints 17,3.5 and 17,5.8
\put(17,3.5){\line(0,1){2.3}}
% endpoints 17,3.5 and 11,2, passes through 12.5,3.5
\put(14.750,-0.250){\arc(2.250,3.750){90}}
% endpoints 11,2 and 9.2,1.5
\put(10.756,-0.613){\arc(0.244,2.613){41.709}}
% Inner arcs
% endpoints 8,2. 8,6.230
\put(8,2){\line(0,1){4.230}}
% endpoints 10,6.5 and 12.5,3.5
\put(13.875,7.188){\arc(-3.875,-0.688){59.490}}
\curvedashes[.1in]{0,1,1}
\put(9,-5.393){\arc(8,7.393){94.518}} % Dashed geodesic
\curvedashes{}
% Arrowheads moved from their theoretical places for better graphics
\put(4.150,4.350){\vector(2,1){0}} % at 4.130,4.350
\put(9.35,5.490){\vector(1,0){0}} % at 9.5.5
\put(12.585,4.891){\vector(3,-1){0}} % at 12.445.4.941
\put(17,2){\vector(1,-1){0}}
% Labels on edges
\put(2.5,2.5){\vclap{\rlap{$g_0$}}}
\put(5,2.5){\vclap{\clap{$g$}}}
\put(8,4){\llap{$g_1^{-1}$}}
\put(5.5,6.1){\vclap{\clap{$g_0^{-1}$}}}
\put(4,6){\vclap{\llap{$g^{-1}$}}}
\put(1.8,5){\vclap{\rlap{$g_1$}}}
\put(8.2,4){\rlap{$g_1$}}
\put(11.2,4){\vclap{\llap{$g_2^{-1}$}}}
\put(11,4.8){\vclap{\rlap{$g_2$}}}
\put(15,4.5){\vclap{\rlap{$g_3^{-1}$}}}
\put(14,3.8){\vclap{\rlap{$g_3$}}}
% Labels inside domains
\put(5.5,4){\vclap{\clap{$\sP$}}}
\put(10,3){\vclap{\clap{$g_1\sP$}}}
\put(13.5,6){\vclap{\rlap{$g_1g_2\sP$}}}
\end{picture}
\caption{The cutting-sequence encoding of a geodesic.  The pictured
geodesic has cutting sequence $\ldots,g_0,g_1,g_2,\ldots$, starting with
$g_1$ from $\sP$.}
\end{figure}
\else %\iffigures
$$
\begin{array}{c}
\hline
~~~~~~~~~~~\mbox{Figure~2.1 about here}~~~~~~~~~ \\
\hline
\end{array}
$$
\fi

The transversality hypothesis says that if $\gamma$ intersects some
translate $g( \sP )$ for $g \in \Gamma$, then $\gamma$
intersects the interior of $g( \sP )$, which
implies that the expansion (\ref{eq22}) is well-defined and unique.
If $\gamma^{-1} = \langle \theta_2 , \theta_1 \rangle$ denotes the
reversed geodesic, then $\bga^{-1} \in G_{\sP,\Gamma}^0$ has the cutting
sequence
\beql{eq24}
C( \bga^{-1} ) = ( \ldots , g'_{-1} , g'_0 , g'_1 , \ldots )
\eeq
in which $g'_i : = (g_{-i} )^{-1}$.

Let $\Sigma_{\sP , \Gamma}^0$ denote the set of all cutting sequences
$C( \bga )$ for $\bga \in G_{\sP,\Gamma}^0$.
This set is invariant under the (forward) shift operator $\sigma ( \{
g_i \} ) = \{ g_{i+1} \}$ since
$$
\sigma( C( \gamma )) = C( g_0 \gamma )~,
$$
and the geodesic $g_0 \gamma$ is in $G_{\sP,\Gamma}^0$.

\begin{defn}
The {\em orbit} or {\em shift-equivalence class} $[C]$ of a cutting
sequence $C$ is the union of all its
forward and backward shifts
\BE
[C] := \bigcup_{k\in\ZZ} \sigma^k(C);
\EE
i.e., it is the smallest shift-invariant set containing $C$.
\end{defn}

Thus $[C(\gamma)]$ is contained in $\Sigma_\sP^0$, but is generally not
a closed set.  The geodesics on the surface $\hpl/\Gamma$ are projections
of geodesics on $\hpl$.  Since the shift operator on cutting sequences
corresponds to a motion of the geodesic by an element of $\Gamma$, the
orbit $[C]$ is an invariant of the projected geodesic.

\begin{defn}
The {\em cutting sequence shift}
$\Sigma_{\sP , \Gamma }$ is the
closure of $\Sigma_{\sP , \Gamma}^0$ in the symbol topology.
That is, $\Sigma_{\sP , \Gamma}$ consists of all symbol sequences
$$
( \ldots , g_{-1} , g_0 , g_1 , \ldots )
$$
such that every finite block $(g_i , g_{i+1} , \ldots , g_{i+k} )$ occurs
in some $C( \gamma )$ for a general position geodesic $\gamma \in
G_{\sP,\Gamma}^0$.
\end{defn}

We now extend the definition of cutting sequences to apply to {\em all}
geodesics $\gamma$ that hit the interior of $\sP$.  The set of cutting
sequences $C( \gamma )$ for a general geodesic
$\gamma=\langle\theta',\theta\rangle$ consists of all symbol sequences
that can be obtained as a limit point (in the sequence topology) of a
sequence of $C( \gamma_j )$ having $\gamma_j \in G_{\sP,\Gamma}^0$ such
that $\gamma_j=\langle\theta'_j,\theta_j\rangle$ with
$\theta'_j\to\theta'$ and $\theta_j\to\theta$ as $j\to\infty$.  There
always exist such convergent sequences $\{\gamma_j : j \geq 1 \}$, because
$\Sigma_{\sP,\Gamma}$ is compact in the symbol topology; hence
$C(\gamma) \neq \emptyset$.  The set $C(\gamma)$ may be infinite for
some geodesics $\gamma$.  All cutting sequences in $C( \gamma )$ are
contained in $\Sigma_{\sP , \Gamma }$; however, the closure operation
used in defining $\Sigma_{\sP,\Gamma}$ allows the possibility that
$\Sigma_{\sP,\Gamma}$ contains some symbol sequences not coming from any
geodesic.

\begin{defn}
The {\em shift-equivalence class} $[C(\gamma)]$ is the union of all orbits of
cutting sequences in $C(\gamma)$.
\end{defn}

We show in the specific case $\Gamma=PSL(2,\ZZ)$ and fundamental domain
$\sF$ that each $[C(\gamma)]$ is a union of a finite number of
orbits of individual cutting sequences (Theorem~\ref{encodingthm}).

The totality of all possible convex polygons $\sP$ that are fundamental
domains of some group $\Gamma$ has an explicit characterization.  A
general treatment appears in Maskit~\cite{Mas71}, which covers non-convex
fundamental domains, and also covers groups of isometries of $\hpl$,
which may include orientation-reversing isometries.
Maskit~\cite[section 2]{Mas71} gives a sufficient condition for $\sP$ to
be a fundamental domain, but for finitely-generated groups where $\hpl
/ \Gamma$ has finite volume and $\sP$ is hyperbolically convex it is a
necessary condition as well.

\subsection{Cutting Sequences for the Modular Surface}
\hsp
We now specialize to cutting sequence expansions for geodesics on the
modular surface $\hpl / PSL (2, \ZZ )$.  The modular group $PSL (2, \ZZ
) = SL (2, \ZZ ) / \pm \hc{I}$ has many different convex fundamental
domains $\sP$ that are hyperbolic quadrilaterals; see Beardon
\cite[Example 9.4.4]{Bea83}.  The {\em standard fundamental
domain} for the modular group acting on $\hpl$ is the hyperbolic
triangle
\beql{eq25}
\sF =  \{ z  :  | z |  \geq 1 ~~~\mbox{and ~~~  $ - \df{1}{2}  \leq  Re (z)
\leq  \df{1}{2} \}$}~,
\eeq
which has one ideal vertex (``cusp'') at $i \In$. When regarded as a
fundamental domain, it is a quadrilateral having an elliptic
vertex added at $z = i$, see Figure 2.2.

\iffigures
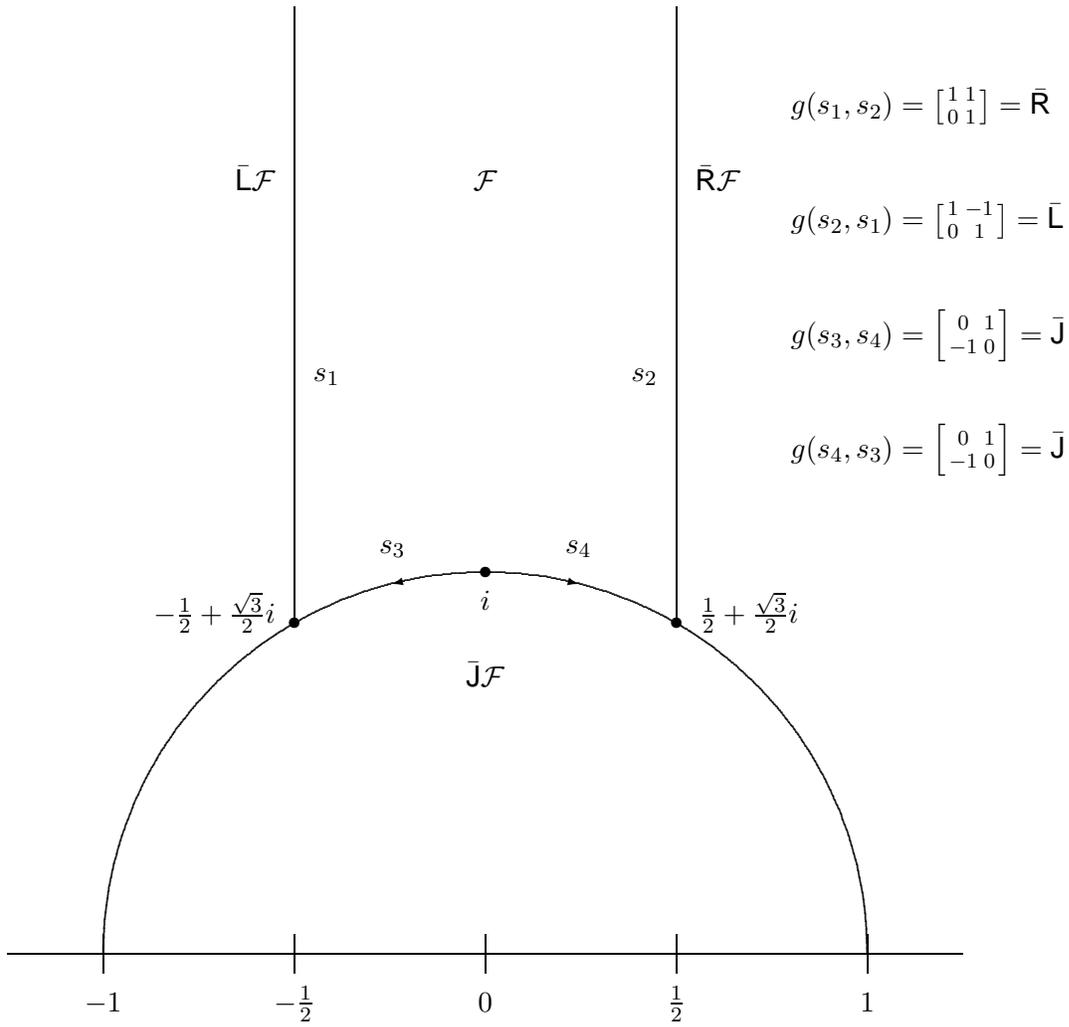
\begin{figure}
%Figure 2.2
\setlength{\unitlength}{2in}
\begin{picture}(2.5,2.5)(-1.25,-.25)
\put(-1.25,0){\line(1,0){2.5}} % x-axis
\multiput(-1,-.05)(.5,0){5}{\line(0,1){.1}} % Tickmarks on x-axis
\put(-1,-.15){\clap{$-1$}} % Labels on x-axis
\put(-.5,-.15){\clap{$-\frac{1}{2}$}}
\put(0,-.15){\clap{0}}
\put(.5,-.15){\clap{$\frac{1}{2}$}}
\put(1,-.15){\clap{1}}
\put(-.5,.866){\line(0,1){1.614}} % vertical lines
\put(.5,.866){\line(0,1){1.614}}
\put(0,0){\arc(1,0){180}} % semicircle
\put(0,2){\clap{$\sF$}} % Label on fundamental domain
\put(-.55,2){\llap{$\bar{\hc{L}}\sF$}} % Labels on adjacent domains
\put(.55,2){\rlap{$\bar{\hc{R}}\sF$}}
\put(0,.7){\clap{$\bar{\hc{J}}\sF$}}
\put(-.45,1.5){\rlap{$s_1$}} % Labels on edges of fundamental domain
\put(.45,1.5){\llap{$s_2$}}
% .243=1/sqrt(17), .970=4/sqrt(17), chosen so that tangent is 1/4
\put(-.243,1.05){\clap{$s_3$}}
\put(.243,1.05){\clap{$s_4$}}
\put(-.243,.970){\vector(-4,-1){0}} % Arrows
\put(.243,.970){\vector(4,-1){0}}
\put(-.5,.866){\circle*{.025}} % Labeled points
\put(-.55,.866){\llap{$-\frac{1}{2}+\frac{\sqrt{3}}{2}i$}}
\put(.5,.866){\circle*{.025}}
\put(.56,.866){\rlap{$\frac{1}{2}+\frac{\sqrt{3}}{2}i$}}
\put(0,1){\circle*{.025}}
\put(0,.9){\clap{$i$}}
% Formulas in top right
\put(.8,2.2){\rlap{$g(s_1,s_2)=\smallmat{1}{1}{0}{1}=\bar{\hc{R}}$}}
\put(.8,1.9){\rlap{$g(s_2,s_1)=\smallmat{1}{-1}{0}{1}=\bar{\hc{L}}$}}
\put(.8,1.6){\rlap{$g(s_3,s_4)=\smallmat{0}{1}{-1}{0}=\bar{\hc{J}}$}}
\put(.8,1.3){\rlap{$g(s_4,s_3)=\smallmat{0}{1}{-1}{0}=\bar{\hc{J}}$}}
\end{picture}
\caption{The fundamental domain $\sF$ with side-pairings.}
\end{figure}
\else %\iffigures
$$
\begin{array}{c}
\hline
~~~~~~~~~~~\mbox{Figure~2.2 about here}~~~~~~~~~ \\
\hline
\end{array}
$$
\fi

The four ``sides'' of $\sF$ are indicated in Figure~2.2 along with the two side
pairings of $s_1$ with $s_2$ and $s_3$ with $s_4$, and
with elements
$g(s, s') \in PSL (2, \ZZ )$, which comprise
\beql{eq26}
N( \sF ) : = N_{PSL(2, \ZZ )}
( \sF ) = \left\{ \left[
\begin{array}{cc}
1 & -1 \\
0 & 1
\end{array}
\right] , \left[
\begin{array}{cc}
1 & 1 \\
0 & 1
\end{array}
\right] , \left[
\begin{array}{cc}
0 & 1 \\
-1 & 0
\end{array}
\right] \right\}~.
\eeq
The ``sides'' $s_3$ and $s_4$ containing the elliptic vertex together make up
one side of $\sF$, and the three
sides are labeled inside the fundamental domain as follows:
\begin{eqnarray*}
\bar{\hc{R}} = \left[ \begin{array}{cc}
1 & 1 \\
0 & 1
\end{array}
\right] & \Leftrightarrow ~~~~~~~s_1 : = \left\{ - \df{1}{2} + it : t >
\df{\sqrt 3}{2} \right\} \\
~~~ \\
\bar{\hc{L}} = \left[ \begin{array}{cc}
1 & -1 \\
0 & 1
\end{array}
\right] & \Leftrightarrow~~~~~~~s_2 : = \left\{ \df{1}{2} + it
: t > \df{\sqrt 3}{2} \right\} ~~~ \\
~~~ \\
\bar{\hc{J}} = \left[
\begin{array}{cc}
0 & 1 \\
-1 & 0
\end{array}
\right] & ~~~~~~~~~~~~~~~~~~~\Leftrightarrow~~~~~~~s_3 \cup s_4 : =
\left\{ |z| = 1 : - \df{1}{2} < \RR e (z) < \df{1}{2} \right\} ~.
\end{eqnarray*}
These matrices
represent the $PSL(2, \ZZ )$-motions needed to move a neighboring fundamental
domain $\sF'$ in $\hpl$ to $\sF$.
These neighboring fundamental domains are labeled with the appropriate
symbol in Figure~2.2.

The action of $PSL (2, \ZZ )$ tiles $\hpl$ with copies of $\sF$, and the
boundaries of all tiles fit together to cover a countable collection of
geodesics of $\hpl$, which are exactly those geodesics having two
rational endpoints $\langle \frac{a}{b}, \frac{c}{d} \rangle$
such that
\beql{eq28}
\det \left[
\begin{array}{cc}
a & c \\
b & d
\end{array}
\right] = \pm 2 ~,
\eeq
where we adopt the convention that the cusp $i \In$ is the rational
$\frac{1}{0}$.
The resulting labeling of the tiles with the labels $\{ \bar{\hc{L}},
\bar{\hc{R}}, \bar{\hc{J}} \}$
is indicated in Figure~2.3.

\iffigures
\begin{figure}
%Figure 2.3
\setlength{\unitlength}{1in}
\begin{picture}(6,4)(-2,-.5)
\put(-2,0){\line(1,0){6}} % x-axis
\multiput(-1,-.1)(1,0){5}{\line(0,1){.2}} % Tickmarks on x-axis
\put(-1,-.2){\clap{$-1$}} % Labels on x-axis
\put(0,-.2){\clap{0}}
\put(1,-.2){\clap{1}}
\put(2,-.2){\clap{2}}
\put(3,-.2){\clap{3}}
\multiput(-1.5,0)(1,0){6}{\line(0,1){4}} %vertical lines
\put(-1,0){\arc(1,0){180}} % semicircles of radius 1
\put(0,0){\arc(1,0){180}}
\put(1,0){\arc(1,0){180}}
\put(2,0){\arc(1,0){180}}
\put(3,0){\arc(1,0){180}}
\put(-2,0){\arc(1,0){90}} % incomplete semicircles of radius 1
\put(4,0){\arc(0,1){90}}
\multiput(-1.667,0)(1,0){6}{\arc(.333,0){180}} % semicircles of radius 1/3
\multiput(-1.333,0)(1,0){6}{\arc(.333,0){180}}
\multiput(-1.6,3)(1,0){6}{\llap{$\bar{\hc{L}}$}} % Labels on vertical edges
\multiput(-1.4,3)(1,0){6}{\rlap{$\bar{\hc{R}}$}}
\put(-1.55,.5){\llap{$\bar{\hcsm{J}}$}}
\put(-1.45,.5){\rlap{$\bar{\hcsm{J}}$}}
\put(-.55,.5){\llap{$\bar{\hcsm{J}}$}}
\put(-.45,.5){\rlap{$\bar{\hcsm{J}}$}}
\put(.45,.5){\llap{$\bar{\hcsm{J}}$}}
\put(.55,.5){\rlap{$\bar{\hcsm{J}}$}}
\put(1.45,.5){\llap{$\bar{\hcsm{J}}$}}
\put(1.55,.5){\rlap{$\bar{\hcsm{J}}$}}
\put(3.45,.5){\llap{$\bar{\hcsm{J}}$}} % No label at x=2.5 because
\put(3.55,.5){\rlap{$\bar{\hcsm{J}}$}} % dashed geodesic in the way
\multiput(-1,1.1)(1,0){5}{\clap{$\bar{\hc{J}}$}} % Labels on bottom of domains
\multiput(-1,.85)(1,0){5}{\clap{$\bar{\hcsm{J}}$}}
\multiput(-1.25,.7)(1,0){6}{\rlap{$\bar{\hcsm{L}}$}} % Labels on semicircles
\multiput(-1.75,.7)(1,0){6}{\llap{$\bar{\hcsm{R}}$}} % of radius 1
\multiput(-1.3,.6)(1,0){6}{\llap{$\bar{\hcsm{R}}$}}
\multiput(-1.7,.6)(1,0){6}{\rlap{$\bar{\hcsm{L}}$}}
\multiput(-1.667,.37)(1,0){6}{\clap{$\bar{\hcsm{R}}$}} % Labels on semicircles
\multiput(-1.333,.37)(1,0){6}{\clap{$\bar{\hcsm{L}}$}} % of radius 1/3
\multiput(-1.667,.22)(1,0){6}{\clap{$\bar{\hcvs{L}}$}}
\multiput(-1.333,.22)(1,0){6}{\clap{$\bar{\hcvs{R}}$}}
\put(0,3.5){\clap{$\sF$}} % Label on fundamental domain
\curvedashes[.1in]{0,1,1}
\put(0,0){\arc(2.5,0){90}} % Dashed geodesic
\curvedashes{}
\put(0,2.5){\circle*{.05}}
\put(0,2.3){\clap{$x$}}
\put(1.768,1.768){\vector(1,-1){0}} % Arrow on geodesic; 1.768=1.25*sqrt(2)
\end{picture}
\caption{Labeled tessellations of $PSL(2,\ZZ)$.  The pictured geodesic,
starting from $x$, has cutting sequence beginning $\bar{\hc{R}},
\bar{\hc{R}}, \bar{\hc{J}}, \bar{\hc{L}}, \bar{\hc{L}}$.}
\end{figure}
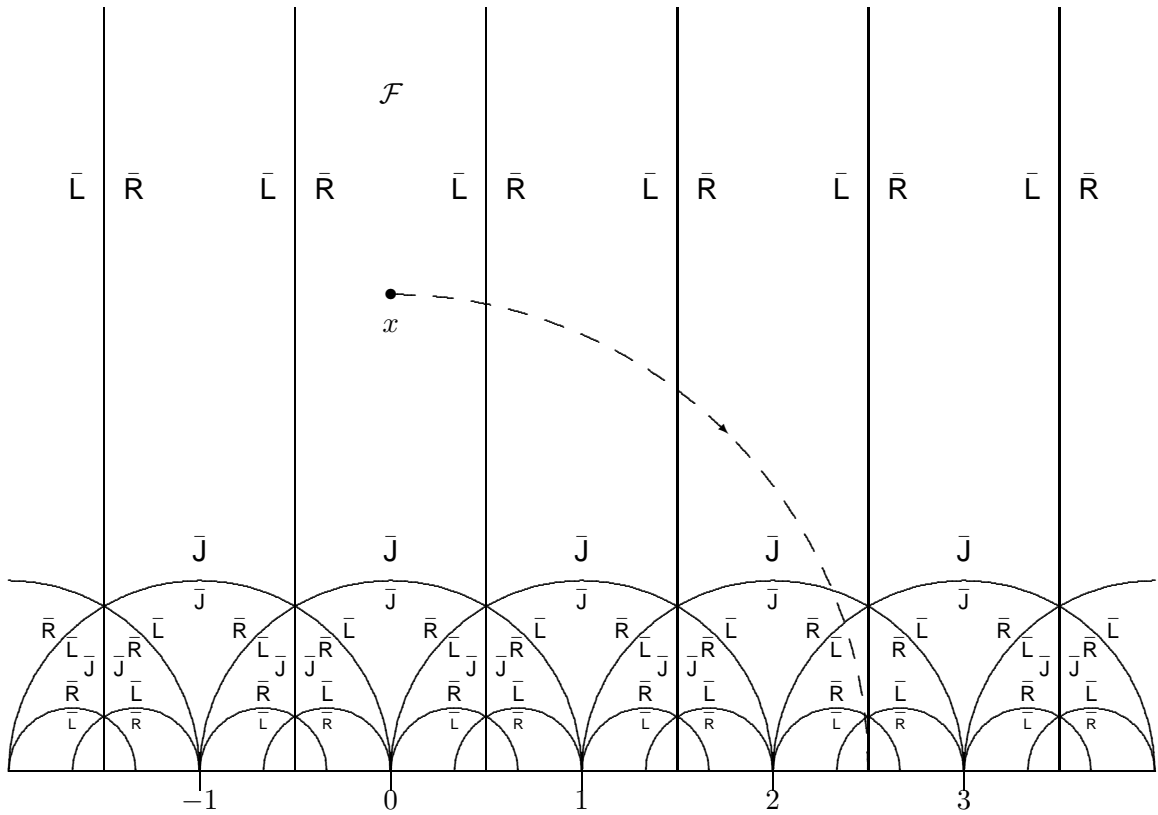
\else %\iffigures
$$
\begin{array}{c}
\hline
~~~~~~~~~~~\mbox{Figure~2.3 about here}~~~~~~~~~ \\
\hline
\end{array}
$$
\fi

The {\em cutting sequence shift} $\Sigma_\sF : = \Sigma_{\sF , PSL (2, \ZZ )}$
is defined by the method of section 2.1,
and is a closed subshift of the shift on three letters $\{ \bar{\hc{L}},
\bar{\hc{R}}, \bar{\hc{J}} \}$.
For this special case, we present some additional cutting sequence labelings
that keep track of corners,
which we call {\em generalized cutting sequences}.
We label the
two finite corners of $\sF$ symbolically by
\begin{eqnarray}
\bar{\hc{C}}_1 = \bar{\hc{J}}\bar{\hc{R}}\bar{\hc{J}} =
\bar{\hc{L}}\bar{\hc{J}}\bar{\hc{L}} =
 \left[ \begin{array}{cc}
-1 & 0 \\
1  &  -1
\end{array} \right] & \Leftrightarrow~~~~~~~\sS_{\hc{C}}^-  : = z =
\df{1}{2} + i \df{\sqrt 3}{2} \\
\bar{\hc{C}}_2 =  \bar{\hc{J}}\bar{\hc{L}}\bar{\hc{J}} =
\bar{\hc{R}}\bar{\hc{J}}\bar{\hc{R}} = \left[ \begin{array}{cc}
-1 & 0 \\
-1 & -1
\end{array}
\right] & ~~~~\Leftrightarrow~~~~~~~\sS_{\hc{C}}^+  : = z = - \df{1}{2} + i
\df{\sqrt 3}{2} ~.
\end{eqnarray}

General geodesics on the modular surface $\hpl / PSL ( 2 , \ZZ )$ are
given by geodesics on $\hpl$ that pass through the fundamental domain.
A representative set of lifts covering every geodesic in the modular
surface consists of all vertical geodesics with $- \frac{1}{2} < \theta
< \frac{1}{2}$, plus all semicircular geodesics in $\hpl$ whose maximum
point lies in the interior of the fundamental domain.  The latter must
necessarily hit the vertical line $\{ it : t > 0 \}$ and have endpoints
$( \theta_1 , \theta_2 )$ satisfying
$$
| \theta_2 - \theta_1 | > \sqrt{3},
\mbox{\ and\ } -1  \leq  \theta_1  +  \theta_2  <1.
$$
These conditions imply that $\theta_1  \theta_2   <  0   $.

A {\em generalized}
{\em cutting sequence}
for a general geodesic $\gamma$
is the sequence of symbols from $\{  \bar{\hc{L}}  ,  \bar{\hc{R}} ,
\bar{\hc{J}} , \bar{\hc{C}}_1 , \bar{\hc{C}}_2   \}$ that
describe the successive sides of the images of the fundamental domains that
it hits.

Finally, for the domain $\sF$ we define a set of {\em vertical cutting
sequences} which are one-sided cutting sequences associated to the set
of oriented
geodesics of $\hpl$ that emanate from the cusp $\{ i \In \}$ of $\sF$
and pass through $\sF$,
which are exactly the set of vertical geodesics
$\{ ( \In , \theta ) : - \frac{1}{2} < \theta < \frac{1}{2} \}$.

\begin{defn}
The {\em irrational vertical cutting sequence set}
$\Pi_\sF^0$ consists of cutting sequences of those vertical geodesics
$\langle\infty,\theta\rangle$ for irrational $\theta$ with
$-\frac{1}{2}<\theta<\frac{1}{2}$.
\end{defn}

An important result for our analysis is the following simple fact.

{\sloppy
\begin{lemma}\label{vertcorner}
An irrational vertical geodesic cannot hit any finite corner
of a $PSL (2, \ZZ )$-translate of $\sF$.
\end{lemma}

}

{\bf Proof.}  The generators of $PSL (2,\ZZ)$ acting on the upper
half-plane are $z\to z+1$ and $z\to -1/z$.  Both generators preserve the
field $\QQ(\sqrt{-3})$, in which every element has rational real part.
The finite corners $\pm1/2 + \sqrt{-3}/2$, and thus all of their images,
are in this field.  $\qed$

Thus the cutting sequences in $\Pi_\sF^0$ are one-sided infinite sequences
which contain no symbol
$\bar{\hc{C}}_1$ or $\bar{\hc{C}}_2$, hence they
are all of the form:
\beql{eq29}
C^+( \theta ) : =
( g_0 , g_1 , g_2 , \cdots ) ~~ \mbox{all ~$g_i \in N( \sF ) =
\{ \bar{\hc{L}}, \bar{\hc{R}}, \bar{\hc{J}} \}$} ~.
\eeq
The expansion $C^+ ( \theta )$ is a kind of continued fraction of $\theta$,
and we call it the {\em cutting sequence expansion} of $\theta$.
Note that, for a given point $\theta + it$ corresponding to a fundamental
domain
$h_i \sF$,
that there is some point $z \in \sF$ such that
\beql{eq2009}
g_0 g_1 \cdots g_{i-2} g_{i-1}(z) = \theta + it,
\eeq
according to Lemma~\ref{translemma}.

\begin{defn}
The {\em vertical cutting
sequence space} $\Pi_\sF$ is the closure of $\Pi_\sF^0$ in the one-sided
sequence topology on the one-sided shift space on three symbols $\{
\bar{\hc{L}}, \bar{\hc{R}}, \bar{\hc{J}} \}$.
\end{defn}

Note that the set $\Pi_\sF^0 $ is not closed under the one-sided shift
operator $\sigma^+$, because the first symbol of every element of
$\Pi_\sF^0$ is $\hc{J}$.

We can now define cutting sequences $C( \theta )$ for rational $\theta$
with $- \frac{1}{2} < \theta < \frac{1}{2}$ in $\Pi_\sF$ as limits of
cutting sequences $C( \theta_i )$ of irrational $\theta_i \rightarrow
\theta$, as described earlier; in this case $C( \theta )$ is a set of
one-sided infinite cutting sequences.

We define {\em generalized cutting sequence expansions} for vertical
geodesics, drawn from the alphabet $\{ \bar{\hc{L}}, \bar{\hc{R}},
\bar{\hc{J}}, \bar{\hc{C}}_1 , \bar{\hc{C}}_2 \}$, similarly to the case
of general geodesics.

\section{Continued Fraction Expansions and Cutting Sequence Expansions}
\hsp
This section describes symbolic dynamics for vertical cutting sequences
associated to additive continued fraction
(ACF) expansions, and for Minkowski geodesic continued fraction expansions
(MGCF) for vertical geodesics.  It shows that the cutting sequence
expansion and
MGCF expansions for a geodesic can each be computed from the other
using a finite automaton.

\subsection{Additive Continued Fraction Expansions}
\hsp
The ordinary continued fraction (OCF) expansion for a real $\theta$, written
$\theta =  [ a_0 , a_1 , a_2 , \cdots ]$,
can be represented using $2 \times 2$ matrices.  The use of such
matrices to represent continued fractions appears in Frame \cite{Fra49}
and Kolden \cite{Kol49}, and is described in Stark \cite{Sta71}.
Namely, one has
\beql{eq201}
\left[
\begin{array}{ll}
a_0 & 1 \\
1 & 0
\end{array}
\right] \ldots \left[
\begin{array}{ll}
a_n & 1 \\
1 & 0
\end{array}
\right] = \left[
\begin{array}{ll}
p_n & p_{n-1} \\
q_n & q_{n-1}
\end{array}
\right],
\eeq
where $\frac{p_n}{q_n} =  [a_0 , a_1  \dd  a_n ]$ is the $n$-th
convergent of $\theta$.
The ordinary continued fraction expansion for $\theta>0$ has a
symbolic dynamics which is a one-sided shift on infinitely many symbols, using
the alphabet
\beql{defineocf}
\hc{L}_k=\twobytwo k 1 1 0 ,\ k\ge 0.
\eeq

The {\em additive continued fraction} (ACF) expansion of $\theta>0$ is a
symbolic expansion using the two symbols
\beql{eq302}
\hc{R} = \left[ \begin{array}{ll}
1 & 1 \\
0 & 1
\end{array} \right] ~,~~
\hc{F} = \left[
\begin{array}{ll}
0 & 1 \\
1 & 0
\end{array}
\right] ~,
\eeq
which  is obtained from~\numb{defineocf} by expanding the symbol
$\hc{L}_k$ as
\BE
\twobytwo k 1 1 0 = {\twobytwo 1 1 0 1}^k \twobytwo 0 1 1 0 .
\EE
The additive continued fraction expansion of $\theta>0$ is thus
\beql{eq202}
\hc{F} \hc{R}^{a_1 } \hc{F} \hc{R}^{a_2 } \hc{F} \hc{R}^{a_3 } \cdots ~,
\eeq
regarded as a string of letters in the alphabet $\{ \hc{R} , \hc{F} \}$.
The finite truncations of this expression are $2 \times 2$ matrices that
encode the convergents and
intermediate convergents of $\theta$.

We next construct a related continued fraction for real $\theta > 0$
which uses the symbols
\beql{eq203}
\hc{R} : = \left[
\begin{array}{ll}
1 & 1 \\
0 & 1
\end{array}
\right]\ , \quad \hc{D} =  \left[
\begin{array}{ll}
1 & 0 \\
1 & 1
\end{array} \right]\ ,
\eeq
The {\em Farey tree expansion}
of a real number $\theta > 0$ that has OCF expansion $\theta = [a_0 ,
a_1 , a_2 , \ldots ]$ is:
\beql{eq204}
\hc{R}^{a_0 } \hc{D}^{a_1 } \hc{R}^{a_2 } \hc{D}^{a_3 }  \cdots
\eeq
regarded as a sequence of letters in the alphabet $\{ \hc{R} , \hc{D} \}$.
This sequence encodes the steps of subtraction needed in the division
process to encode the ordinary continued fraction of $\theta$; it also
essentially gives all the
intermediate convergents to $\theta$, cf. Richards \cite{Ric81}, Theorem 2.1.
For $\theta > 0$ the set of allowable symbol sequences for the additive
expansion
is the full one-sided shift on two letters $\{  \hc{R} , \hc{D}  \}$.

A relation of this expansion to paths in the Farey tree is described in
Lagarias~\cite{Lag92} and Lagarias and Tresser~\cite{LagTre95}.
We start with the interval
$\bigl[\frac{0}{1},\frac{1}{0}]$, represented as the matrix $\twobytwo 0
1 1 0$.  At each step, we multiply our current matrix $\twobytwo
{p_1}{q_1}{p_2}{q_2}$ on the left by $\hc{R}$ or $\hc{D}$, as
appropriate.  The matrix $\hc{R}$ causes our current matrix $\twobytwo
{p_1}{q_1}{p_2}{q_2}$ to be replaced by
$\twobytwo{p_1+p_2}{q_1+q_2}{p_2}{q_2}$; the matrix $\hc{D}$ causes it
to be replaced by $\twobytwo{p_1}{q_1}{p_1+p_2}{q_1+q_2}$,  Thus, in
either case, the next term in the Farey tree is $(p_1+p_2)/(q_1+q_2)$,
and we continue with the interval on the side of this approximation
which contains $\theta$.

The symbolic expansions (\ref{eq202}) and (\ref{eq204}) are equivalent in
the sense that each is convertible
to the other using a finite automaton; see section~\ref{automata}.

\subsection{Cutting Sequences and Minkowski lattice basis reduction}
We now develop a correspondence between the cutting sequence expansion
of the geodesic $\theta+it$ for $- \frac{1}{2} <  \theta  < \frac{1}{2}$
and the construction of a Minkowski-reduced lattice basis for the
parametrized series of lattice bases
\beql{eq205}
\hc{B}_t ( \theta ) =  \left[
\begin{array}{cc}
1 & 0 \\
- \theta & t
\end{array}
\right]
\eeq
for the parametrized family of lattices $\Lambda_t =  \ZZ [ ( 1,0), (-
\theta , t)]$.

To demonstrate the equivalence, we first transform lattice bases in $GL
(2, \RR )$ to the cone of positive definite symmetric matrices $\sP_2$
(equivalently, to positive definite quadratic forms), and then we relate
these to elements of $\hpl$.

\begin{defn}
A lattice basis $\bv_1,\bv_2$ is {\em Minkowski-reduced} if $\bv_1$ is a
shortest vector in the lattice, and $\bv_2$ is a shortest vector not a
multiple of $\bv_1$.
\end{defn}

For each $t$, there is a matrix $\hc{P}=\hc{P}_t$ (possibly several) for
which $\hc{P}\hc{B}_t(\theta)$ is Minkowski-reduced.
General properties  of Minkowski-reduction can be found in Cassels~\cite{Cas78}
and Gruber and Lekkerkerker~\cite[pp. 149ff]{GL87}.

A basis $\hc{B} \in GL(2, \RR )$ corresponds to the
positive definite symmetric matrix $\hc{M} = \hc{B} \hc{B}^T \in \sP_2$
with associated quadratic form
\beql{eq3014}
Q(x_1 , x_2 ) = [x_1 x_2 ] \hc{M} \left[ \begin{array}{c}
x_1 \\
x_2
\end{array}
\right] = m_{11} x_1^2 + (m_{12} + m_{21} )
x_1 x_2 + m_{22} x_2^2 ~.
\eeq
For the basis matrix $\hc{B}_t ( \btheta )$ in (\ref{eq205}) this gives
\beql{eq3015}
\hc{M} : = \hc{M}_t ( \theta ) = \left[ \begin{array}{cc}
1 & - \theta \\
- \theta & \theta^2 + t^2
\end{array} \right]
\eeq
with associated quadratic form
\beql{eq3016}
Q_t (x_1 , x_2 ) = ( x_1 - \theta x_2 )^2 +
t^2 x_2^2 = (x_1 - ( \theta+ it) y_1 ) (x_1 -
( \theta - it ) x_2 ) ~.
\eeq
The action of $\hc{P} \in GL(2, \ZZ )$ which sends the matrix
$\hc{B}_t ( \theta ) $ to $\hc{P}\hc{B}_t ( \btheta )$
is transformed to the $GL(2, \ZZ )$-action on $\sP_2$ that takes $\hc{M}$ to
$\hc{P} \hc{M} \hc{P}^T$, and this takes the associated quadratic
form $Q(x_1 , x_2 )$ to
\beql{eq3017}
Q' (x_1 , x_2 ) = Q( p_{11} x_1 + p_{21} x_2 , ~
p_{21} x_1 + p_{22} x_2 ) ~.
\eeq
The Minkowski reduction domain $\bar{\sM}$ for $GL( 2, \ZZ ) \backslash \sP_2$
(see Cassels \cite{Cas78}) is the set of quadratic forms~\numb{eq3014}
such that
\begin{eqnarray}
Q(1,0) & \leq & Q( 0,1) ~, \nonumber \\
Q(0,1) & \leq & Q(1,1) ~, \label{quadineqs}\\
Q(0,1) & \leq & Q(-1 , 1) ~. \nonumber
\end{eqnarray}
This is equivalent to the conditions
\beql{eq3019}
|m_{12} + m_{21} | \leq m_{11} \leq m_{22} ~.
\eeq
on the form $Q$ in (\ref{eq3014}).
The subgroup
$$
H_2 : = \left\{ \left[
\begin{array}{cc}
1 & 0 \\
0 & 1
\end{array} \right] ~, ~
\left[
\begin{array}{cc}
-1 & 0 \\
0 & 1
\end{array}
\right] , ~ \left[
\begin{array}{cc}
1 & 0 \\
0 & -1
\end{array} \right] ~,~ \left[ \begin{array}{cc}
-1 & 0 \\
0 & -1
\end{array}
\right] \right\}
$$
of $GL(2, \ZZ )$ maps the domain $\bar{\sM}$ into itself,
hence we actually let $H_2 \backslash GL(2, \ZZ )$ act on $\sP_2$.
By definition the Minkowski reduction domain $\sM$ of
$GL(2, \ZZ ) \backslash GL(2, \RR)$,
is the preimage of $\bar{\sM}$ under the mapping $\bar{\hc{B}} \mapsto
\hc{B} \hc{B}^T$,
and the group actions by $GL(2, \ZZ )$ coincide.
An algorithm for choosing Minkowski-reduced lattice bases will choose
one element of $GL(2, \ZZ )$ out of the four-element coset of $H_2
\backslash GL(2, \ZZ )$.  A natural canonical choice is to
require $\det \hc{P}=1$; this still leaves two choices for $\hc{P}$, but
multiplying $\hc{P}$ by $-1$ will not change the quadratic form and can
be done arbitrarily.  A different choice will be made by the Minkowski
geodesic continued fraction of the next section.

For the second transformation from $\sP$ to $\hpl$, note that a
positive definite quadratic form can be uniquely written
\beql{eq320}
Q(x,y) = a(x_1 - \theta x_2 ) (x_1 - \bar{\theta} x_2 ) ~,
\eeq
with $\theta$ and $\bar{\theta}$ complex conjugates and $Im ( \theta ) > 0$.
We map $\sP$ to $\hpl$ by sending the quadratic
form $Q$ to the unordered pair of roots $\{ \theta, \bar{\theta} \}$ and
then identify this with
$\theta \in \hpl = \{ z : Im (z) > 0 \}$.
If $\hc{P} \in GL(2, \ZZ )$ sends $Q$ to $Q'$ given by~\numb{eq3019},
then a calculation yields
$$
Q' = a'   (x_1 - \phi x_2 ) (x_1 - \bar{\phi} x_2 ) ~,
$$
with $a' = Q(p_{11} , p_{12})$ and with a root $\phi$ given by
\beql{eq312}
\phi = ( \hc{P}^T )^{-1} \theta = \left[
\begin{array}{cc}
p_{22} & - p_{21} \\
-p_{12} & p_{11}
\end{array}
\right] ( \theta ) ~,
\eeq
where the action of the matrix is a linear fractional transformation
on $\CC$.
If $\det ( \hc{P} ) = 1$, then $\phi \in \hpl$, while if
$\det ( \hc{P} ) = -1$ then
$Im ( \phi ) < 0$ is the complex conjugate of the root we want, i.e.
it reverses the ordering of the roots.
Now each four element coset of $H_2 \backslash GL(2, \ZZ )$ contains
two elements in $SL(2, \ZZ )$ which form a coset of
$PSL (2, \ZZ ) = \{ \pm \hc{I} \} \backslash SL(2, \ZZ )$, so the
action~\numb{eq312} of $H_2 \backslash GL(2, \ZZ )$ on unordered pairs $\{
\theta , \bar{\theta} \}$
is equivalent to the standard $PSL (2, \ZZ )$-action on $\hpl$.
Finally, the conditions~\numb{quadineqs} for a form $Q$ to be Minkowski
reduced translate via~\numb{eq320} to
\[
\theta \bar{\theta} \geq 1,\quad  -1 \geq \theta + \bar{\theta} \geq 1 ~,
\]
which is exactly the fundamental domain $\sF$ of $PSL (2, \ZZ )$.

The mapping from $GL(2, \RR )$ to $\hpl$ obtained by composing these
two transformations sends the set of matrices $\{ \hc{B}_t ( \theta ) :
t > 0 \}$ in
$GL(2, \RR )$ to the geodesic $\{ \theta + it : t > 0 \}$ in $\hpl$.
Since the $H_2 \backslash GL(2, \ZZ )$ action on $GL( 2, \RR )$ is equivalent
to the $PSL (2, \ZZ )$ action on $\hpl$ and the reduction domains correspond,
the sequence of Minkowski-reduced lattice bases for $\theta$ is
essentially the same as the cutting sequence expansion for $\theta$.
In the reverse direction, $\theta \in \hpl$ determines a positive ray
$\{ aQ : a > 0 \}$ in the cone $\sP_2$, i.e. an element of $\sP_2 / \RR^+$,
and this in term determines a family $\{ a B Q : a \in \RR^+ , Q \in O (n,
\RR ) \}$
in $GL(2, \RR )$, i.e. an element of $GL(2, \RR )/ \RR^+ O(2, \RR )$.
This does not affect the symbolic dynamics because the Minkowski domain
$\bar{\sM}$ is
invariant under the $\RR^+$-action, and because the Minkowski reduction domain
$\sM$ is invariant under the
$R^* O (2, \RR )$-action.

\subsection{Minkowski Geodesic Continued Fraction Expansions}\label{secmink}
\hsp
The Minkowski geodesic multidimensional continued fraction (MGCF)
expansion is introduced in Lagarias~\cite{Lag94}.
We consider here the one-dimensional case.

This is a specific algorithm for Minkowski-reduction of the lattice bases
\BE
\hc{B}_t ( \theta ) =  \left[
\begin{array}{cc}
1 & 0 \\
- \theta & t
\end{array}
\right]
\EE
There exists a sequence of
{\em critical values}
$$
\In =  t_0  >  t_1 >  t_2  >  \cdots
$$
with $\lim_{n \rightarrow \In} t_n =  0$, and an
associated sequence of
{\em convergent matrices}
$\{  \hc{P}^{(n)}  :  n =  0,1,2, \cdots \} $ in $GL( 2, \ZZ )$ such that
$\hc{P}^{(n)}  \hc{B}_t ( \theta )$ is a Minkowski-reduced lattice basis of the
lattice when $t_n  >  t  >  t_{n+1}$.
We write
\beql{eq36}
 \hc{P}^{(n)}:  = \left[
\begin{array}{ll}
p_1^{(n)} & q_1^{(n)} \\
p_2^{(n)} & q_2^{(n)}
\end{array}
\right] ~,
\eeq
so that
\beql{eq206}
\hc{P}^{(n)} \hc{B}_t ( \theta ) =  \left[
\begin{array}{ll}
p_1^{(n)} - q_1^{(n)} \theta & q_1^{(n)} t \\
p_2^{(n)} - q_2^{(n)} \theta & q_2^{(n)} t
\end{array}
\right] ~.
\eeq
This motivates the name ``continued fraction''; as $t$ goes to zero, the
approximations $p_i^{(n)} - q_i^{(n)}\theta$ must also go to zero, and
thus the $p_i/q_i$ play the role of convergents.

Changing the sign of one or both basis vectors does not affect
Minkowski-reduction, so we have four choices of $\hc P^{(n)}$.  For
defining the MGCF as a continued fraction, it is natural to require
positive denominators; that is, we require that $q_i^{(n)} \geq 0$ for
all $n \geq 1$, and if $q_i^{(n)} = 0$ then $p_i^{(n)} > 0$.  For
lattice basis reduction, it is natural to require $\det \hc P^{(n)}=1$;
multiplying by $\pm\hc I$ is still possible but will not affect the
process.

In either case, the associated {\em partial quotient matrices}
$\hc{A}^{(n)}$ are defined by
$$
\hc{P}^{(n)} =  \hc{A}^{(n)} \hc{P}^{(n-1)} ~,
$$
so that
\beql{eq207}
\hc{P}^{(n)} =  \hc{A}^{(n)} \hc{A}^{(n-1)} \cdots \hc{A}^{(1)} \hc{P}^{(0)} ~.
\eeq
Here $\hc{P}^{(0)}  = \hc{I}$ is the identity matrix if $- \frac{1}{2}
\leq  \theta  <  \frac{1}{2}$, since the lattice basis
$\{(1,0),(-\theta,t)\}$ is Minkowski-reduced for large $t$.

We now show that the partial quotient matrices $\hc{A}^{(n)}$ are drawn
from a finite set.

\begin{lemma}[Minkowski partial quotient set].
The allowed partial quotients for the one-dimensional Minkowski geodesic
continued fraction algorithm are
\beql{eq310}
\left\{ \left[
\begin{array}{cc}
0 & 1 \\
1 & 0
\end{array}
\right] , ~ \left[
\begin{array}{cc}
1 & 0 \\
1 & 1
\end{array}
\right] , ~ \left[
\begin{array}{cc}
1 & 0 \\
1 & -1
\end{array}
\right] , ~ \left[
\begin{array}{cc}
1 & 1 \\
0 & 1
\end{array}
\right] \right\} ~.
\eeq
\end{lemma}

{\bf Proof.}
Let
$\hc{M} = \left[ \begin{array}{c}
\bv_1 \\
\bv_2
\end{array}
\right]$
be a $2 \times 2$ matrix
with $\det (  \hc{M}  )  >  0$, and
associate to it the positive definite quadratic form with coefficient matrix
\beql{eq209}
 \hc{M}   \hc{M}^T =  \left[
\begin{array}{ll}
\| \bv_1 \|^2 & \langle \bv_1 , \bv_2 \rangle \\
\langle \bv_1 , \bv_2 \rangle & \| \bv_2 \|^2
\end{array}
\right] ~.
\eeq
The Minkowski-reduction conditions for the quadratic form $Q( \bx ) =
\bx^T (  \hc{M}   \hc{M}  ^T ) \bx$ are
\begin{eqnarray}
\| \bv_1  \| & \leq & \| \bv_2 \| \label{mink1}\\
\| \bv_2 \|  & \leq & \| \bv_1   +  \bv_2 \| \label{mink2}\\
\| \bv_2 \|  & \leq  & \| \bv_1   - \bv_2  \| ,\label{mink3}
\end{eqnarray}
see Cassels \cite[p. 257]{Cas78}.

We must now choose our algorithm; we will describe both the natural
algorithm for Min\-kowski lattice basis reduction and the Minkowski
geodesic continued fraction algorithm.

For lattice basis reduction, it is natural to apply a new matrix
$\hc{A}^{(n)}$ of determinant 1.  At any critical time $t_n$, all of the
above inequalities are satisfied with the current matrix $\hc{P}^{(n-1)}$
for $t_{n-1}>t>t_n$, and at least one holds with equality at $t=t_n$,

We first consider ``generic'' convergents which occur when exactly two
of $\| \bv_1 (t) \|$,
$\| \bv_2 (t) \|$,
$\| \bv_1 (t) + \bv_2 (t) \|$ and
$\| \bv_1 (t) - \bv_2 (t) \|$ become equal at
$t = t_n$.

If at $t = t_n$ we have $\| \bv_1 (t ) \| = \| \bv_2 (t ) \|$ then the
partial quotient matrix is
$\tilde{\hc{J}}=\left[ \begin{array}{cc}
0 & 1 \\
-1 & 0
\end{array}
\right]$; it exchanges $\bv_1$ and $\bv_2$ and changes the sign of the
new $\bv_2$ to keep determinant 1.
If $\| \bv_2 (t) \| = \| \bv_1 (t) + \bv_2 (t) \|$ then the partial
quotient matrix is $\tilde{\hc{R}}=\left[
\begin{array}{cc}
1 & 0 \\
1 & 1
\end{array}
\right]$;
it replaces $\bv_2$ with $\bv_1 + \bv_2$.
If $\| \bv_2 (t) \| = \| \bv_1 (t) - \bv_2 (t) \|$ then the partial
quotient matrix is $\tilde{\hc{L}}=\left[
\begin{array}{cc}
1 & 0 \\
-1 & 1
\end{array}
\right]$;
it replaces $\bv_2$ with $\bv_2 - \bv_1$.

We might also have more than one of the
inequalities~\numb{mink1}--\numb{mink3} holding with equality at the
same $t$.  If $\|\bv_1(t)+\bv_2(t)\|=\|\bv_1(t)-\bv_2(t)\|$, then neither one
can be equal to $\|\bv_1(t)\|$ unless $\bv_2(t)=0$, nor to $\|\bv_2(t)\|$
unless $\bv_1(t)=0$.  We could have $
\| \bv_1(t) \| = \| \bv_2(t) \| = \| \bv_1(t) + \bv_2(t)  \|$.  In that case, for
$t<t_n$, we would have $
\| \bv_1 (t) \| > \| \bv_2 (t) \| > \| \bv_1 (t) + \bv_2 (t) \|$, and to
correct these inequalities, we must replace $\bv_1$ by $\bv_1+\bv_2$;
we can take partial quotient matrix
$\tilde{\hc{C}}_1=\tilde{\hc{J}}\tilde{\hc{L}}\tilde{\hc{J}}=\left[
\begin{array}{cc}
-1 & -1 \\
0 & -1
\end{array}
\right]$ by also changing signs of both vectors   Similarly, if $
\| \bv_1(t) \| = \| \bv_2(t) \| = \| \bv_1 (t) - \bv_2 (t) \|$, we get
partial quotient matrix
$\tilde{\hc{C}}_2=\tilde{\hc{J}}\tilde{\hc{R}}\tilde{\hc{J}}=\left[
\begin{array}{cc}
-1 & 1 \\
0 & -1
\end{array}
\right].$

These five matrices (generated from just three) give the
natural algorithm for Minkowski lattice basis reduction.

For the Minkowski geodesic continued fraction, we follow a similar
process except that we choose signs to keep the $q_i$ positive.
We use the fact that all vectors in the lattice $\Lambda_t ( \theta )$ have
the special form
\BE
\bv_i = \bv_i (t) : = (p_i - q_i \theta , q_i t) {\rm\  with\ }
q_i \geq 0.
\EE
Whenever a critical value $t_n$ occurs with
$
\| (p-q \theta , qt_n ) \| = \| ( p'-q' \theta , q't_c ) \|
$,
then the inequality
\[
\| (p'-q' \theta , q't ) \| > \| ( p-q \theta , qt ) \| ~~\mbox{for ~$0 < t
< t_c$}
\]
holds exactly when $|q'| > |q|$.
Thus when the shortest vector in the basis $\hc{P}^{(n-1)}$ is replaced
by a vector in $\hc{P}^{(n)}$ the associated denominator must increase.

Again, we first consider ``generic'' convergents.

Suppose first that $q_1 < q_2$.
If at $t = t_n$ we have $\| \bv_1 (t ) \| = \| \bv_2 (t ) \|$ then the
partial quotient is
$\left[ \begin{array}{cc}
0 & 1 \\
1 & 0
\end{array}
\right]$; it exchanges $\bv_1$ and $\bv_2$.
If $\| \bv_2 (t) \| = \| \bv_1 (t) + \bv_2 (t) \|$ then the partial
quotient is $\left[
\begin{array}{cc}
1 & 0 \\
1 & 1
\end{array}
\right]$;
it replaces $\bv_2$ with $\bv_1 + \bv_2$.
The case $\| \bv_2 (t) \| = \| \bv_1 (t) - \bv_2 (t) \|$
cannot occur.

Suppose next that $q_1 \geq q_2$.
At $t = t_n$ we cannot have $\|\bv_1 (t) \| = \| \bv_2 (t) \|$ since an
exchange of $\bv_1$ and $\bv_2$ would not increase the denominator.
We may have $\|\bv_1 (t) + \bv_2 (t) \| = \| \bv_2 (t) \|$, giving the
partial quotient $\left[
\begin{array}{cc}
1 & 0 \\
1 & 1
\end{array}
\right]$.
If $q_1 > 2q_2$ then $q' = q_1 - q_2 > q_2$ hence the
partial quotient $\left[
\begin{array}{cc}
1 & 0 \\
1 & -1
\end{array}
\right]$
is possible; it replaces $\bv_2$ with $\bv_1 - \bv_2$.

There remain situations in which three or more of
$\| \bv_1 (t) \|$, $\|\bv_2 (t) \|$,
$\| \bv_1 (t) + \bv_2 (t) \|$ and
$\| \bv_1 (t) - \bv_2 (t) \|$ simultaneously becoming equal
at $t = t_c$.
Only one case is possible; it is
$$
\| \bv_1 (t) \| = \| \bv_2 (t) \| = \| \bv_1 (t) + \bv_2 (t) \| ~,
$$
which can only occur when $q_1 < q_2$.
In this case the partial quotient is $\left[
\begin{array}{cc}
1 & 1 \\
0 & 1
\end{array}
\right]$.  The other case consistent with an increasing denominator is
$\| \bv_2(t) \| = \| \bv_1(t) + \bv_2(t) \| = \| \bv_1(t) - \bv_2(t) \|$,
but this implies $\| \bv_1(t)\| = 0 $, which is impossible unless $t=0$.
$\qed$

We define the
{\em Minkowski geodesic continued fraction expansion}
of any real $\theta$ satisfying
$- \frac{1}{2}  <  \theta  <  \frac{1}{2}$ to be the symbol sequence
\beql{eq208}
( \hc{A}^{(1)}, \hc{A}^{(2)}, \hc{A}^{(3)} , \ldots ) ~.
\eeq
We write this sequence left-to-right, although the matrix product
in~\numb{eq207} runs right-to-left.
The associated {\em Minkowski geodesic symbol set}
is
\beql{eq313}
\hc{L} = \left[  \begin{array}{cc}
1 & 0 \\
1 & -1
\end{array} \right] , \hc{R} = \left[
\begin{array}{cc}
1 & 0 \\
1 & 1
\end{array} \right] , \hc{J} = \left[
\begin{array}{cc}
0 & 1 \\
1 & 0
\end{array} \right] , \hc{C} = \left[
\begin{array}{cc}
1 & 1 \\
0 & 1
\end{array}
\right] ~.
\eeq
The symbol $\hc{C}$ occurs only for ``exceptional'' geodesics, and we will
usually be
concerned with symbolic expansions drawn from the smaller symbol set
$\{ \hc{L} , \hc{R}, \hc{J} \}$.

The Minkowski geodesic continued fraction expansion corresponding to a
general geodesic on $\sF$
is the expansion attached to Minkowski reduction of the parametrized
lattice bases
$$
\hc{B}_t  ( \theta_1 , \theta_2 ) =  \left[
\begin{array}{ll}
1 & (\theta_1 )^{-1} t \\
- \theta_2 & t
\end{array}
\right] ~.
$$
The associated quadratic form is
\beql{eq314}
Q_t ( x,y)  =  (x - \theta_2 y )^2  +  t^2 \left( \df{1}{\theta_1 } x-y
\right)^2 ~,
\eeq
and the partial quotients are obtained by the same formula.

\subsection{Correspondence between the continued fraction and cutting sequence}

We will now give a correspondence between the symbol
sequences for $\theta$ given by cutting sequences, the natural
Minkowski lattice basis reduction algorithm, and the Minkowski geodesic
continued fraction algorithm.

The precise correspondence between the symbol sequences involves
specifying the relation among the four elements of the $H_2$-coset in
$H_2 \backslash GL(2, \ZZ )$ for the Minkowski geodesic continued
fraction expansion and the two elements in the $\{ \pm \hc{I} \}$-coset
in $PSL (2, \ZZ )$ for the cutting sequence expansion.

\begin{defn}
The {\em parity} of a word $\hc{W}=\hc{S}_1\ldots\hc{S}_n$ in the
alphabet $\{\hc{L},\hc{R},\hc{J}\}$ is even or odd according to whether
there are an even or odd number of $\hc{L}$ and $\hc{J}$.  That is,
\BE
\det(\hc{W})=(-1)^{\textstyle\rm parity(\hc{W})}.
\EE
\end{defn}

We obtain:

\begin{thm}\label{MGCFtoCS}
For irrational $\theta$ with $- \frac{1}{2} < \theta < \frac{1}{2}$, the
one-sided cutting sequence expansion of the vertical geodesic $\langle
\infty , \theta \rangle = \{ \theta + it : t > 0 \}$ in $\Pi_\sF^0$ is
obtained from the Minkowski geodesic continued fraction expansion of
$\theta$ by the following procedure: If the current initial word of the
MGCF expansion has even parity, on the next symbol make the letter
replacement $\hc{L} \rightarrow \bar{\hc{L}} , \hc{R} \rightarrow
\bar{\hc{R}}$ and $\hc{J} \rightarrow \bar{\hc{J}}$; if it has odd
parity, make the letter replacement $\hc{L} \rightarrow \bar{\hc{R}} ,
\hc{R} \rightarrow \bar{\hc{L}}$ and $\hc{J} \rightarrow \bar{\hc{J}}$.
The MGCF can be obtained from the cutting sequence by the reverse
process, in which the parity is determined by the current symbols of the
MGCF expansion.
\end{thm}

{\bf Proof.}
By Lemma~\ref{vertcorner}, the assumption of irrational $\theta$ ensures
that the MGCF expansions and cutting sequence expansions of $\theta$ are
both infinite and never use a symbol
$\hc{C}$.  The MGCF expansion
$\cdots\hc{A}^{(3)}\hc{A}^{(2)}\hc{A}^{(1)}$ uses the symbols
$$
\hc{L} = \left[
\begin{array}{cc}
1 & 0 \\
1 & -1
\end{array} \right] , \hc{R} = \left[
\begin{array}{cc}
1 & 0 \\
1 & 1
\end{array} \right] ~,~
\hc{J} = \left[
\begin{array}{cc}
0 & 1 \\
1 & 0
\end{array}
\right] ~,
$$
drawn from $GL(2, \ZZ )$ and moves right to left, while the cutting
sequence expansion uses the symbols
$$
\bar{\hc{L}} = \left[
\begin{array}{cc}
1 & -1 \\
0 & 1
\end{array}
\right] ~,~~ \bar{\hc{R}} = \left[
\begin{array}{cc}
1 & 1 \\
0 & 1
\end{array}
\right] ~ , ~ \bar{\hc{J}} = \left[
\begin{array}{cc}
0 & 1 \\
-1 & 0
\end{array}
\right]
$$
drawn from $SL(2, \ZZ )$ and moves left to right by Lemma~\ref{translemma}.
For notational convenience, we introduce the matrix
\[
\hc{K} := \left[
\begin{array}{cc}
1 & 0 \\
0 & -1
\end{array}
\right],
\]
so $H_2$ is $\{\hc{I},\hc{K},-\hc{K},-\hc{I}\}$.
We first arrange for matrices to multiply in the same direction.  The
matrices in the MGCF are chosen so that
$\hc{P}^{(n)}B_t(\theta)\in\sF$ for a given $t$, while the matrices in
the cutting sequence are chosen so that the cutting sequence product
$h_n=\hc{S}_1\ldots\hc{S}_n$ has $\gamma(t)\in h_n\sF$.  Therefore, we
have $h_n = ( \hc{P}^{(n)} )^{-1}$.
We expand $( \hc{P}^{(n)} )^{-1}$ in terms of symbols
$\hc{L}, \hc{R}, \hc{J}, \hc{K}$ using the
relations $\hc{L}^{-1} = \hc{L}$,
$\hc{R}^{-1} = \hc{K} \hc{R} \hc{K}$ and
$\hc{J}^{-1} = \hc{J}$ to obtain an expansion which proceeds
left to right.  To convert this expansion to the cutting sequence
expansion, we must
convert to symbols $\bar{\hc{L}}, \bar{\hc{R}}, \bar{\hc{J}}$ and remove
the symbols $\hc{K}$, which encode the parity.  We first convert to the
symbols of the natural Minkowski basis reduction algorithm of
Section~\ref{secmink}.
\[
\tilde{\hc{R}} : = \hc{R},\quad
\tilde{\hc{L}} : = \left[ \begin{array}{cc}
1 & 0 \\
-1 & 1
\end{array}
\right], \quad
\tilde{\hc{J}} : = \left[
\begin{array}{cc}
0 & 1 \\
-1 & 0
\end{array}
\right]
\]
in $SL(2,\ZZ)$.  To convert a product of matrices from the form
$\hc{L},\hc{R},\hc{J}$ to
$\tilde{\hc{L}},\tilde{\hc{R}},\tilde{\hc{J}}$, starting from the right
end of the product, we use the relations $\hc{L}\hc{K}=\tilde{\hc{R}}$,
$\hc{R}\hc{K}=\hc{K}\tilde{\hc{L}}$, and
$\hc{J}\hc{K}=-\tilde{\hc{J}}$.
In doing this we pick up or lose a factor of $\hc{K}$ whenever we
encounter a matrix $\hc{L}$ or $\hc{J}$ and this multiplies the
determinant by $-1$; meanwhile, every $\hc{R}$ in the MGCF becomes
$\hc{R}^{-1}=\hc{K}\hc{R}\hc{K}$ when inverted, and thus is itself
encoded with the opposite parity but leaves the parity unchanged for the
next matrix.  (We sometimes pick
up a matrix factor of $-\hc{I}$, but this commutes with everything and
may be  ignored.)  Similarly, to reverse this, we use
$\tilde{\hc{L}}=\hc{K}\hc{L}$ and $\tilde{\hc{J}}=\hc{K}\hc{J}$.

To convert from the natural Minkowski basis reduction algorithm symbols
$\tilde{\hc{L}},\tilde{\hc{R}},\tilde{\hc{J}}$ to the cutting sequence
symbols $\bar{\hc{L}},\bar{\hc{R}},\bar{\hc{J}}$, we use the relations $(
\tilde{\hc{L}}^T )^{-1} = \bar{\hc{R}}$, $(
\tilde{\hc{R}}^T )^{-1} = \bar{\hc{L}}$ and $( \tilde{\hc{J}}^T )^{-1} =
\bar{\hc{J}}$.  Thus the conversion from
$\tilde{\hc{L}},\tilde{\hc{R}},\tilde{\hc{J}}$ to
$\bar{\hc{L}},\bar{\hc{R}},\bar{\hc{J}}$ interchanges $\hc{R}$ with
$\hc{L}$.  $\qed$

\subsection{Finite Automata}\label{automata}

By a {\em finite automaton} we mean a deterministic finite-state automaton, as
defined in Hopcroft and Ullman \cite{HopUll79} or \cite{Ran73},
used as a transducer. Such an automaton is a finite directed graph
with labeled edges, which may contain loops and several edges exiting
from each vertex. The states are the vertices of the graph,
and the edges give rules to move from one state
to the next. Each edge has two labels, an input label and an
output label. If the symbol alphabet has
$s$ letters, then from each vertex there exit exactly $s$ edges whose
input labels are exactly the $s$ allowed symbols. The output labels
are finite strings of symbols, possibly empty. The machine starts
in a given state. It reads an input symbol, which tells it which exit
edge to follow, prints the specified output string of letters as 
output, and moves to the new state specified by the edge.
Then it proceeds  to the next input symbol.

In this paper we present a number of results asserting the existence of
finite automata to convert one-sided infinite symbol sequences of one form to
another form.
In the proofs we only indicate the ``finite-state'' character of the
conversion process,
and generally omit details of the (routine but sometimes involved)
construction of the
automaton.

As a simple example, the discussion in section 3.1 yields:
\begin{thm}\label{ACFtoFarey}
For real $\theta>1$, the additive continued fraction expansion
of $\theta$ can be converted to the Farey tree expansion of $\theta$ by
a finite automaton, and vice versa.
\end{thm}

{\bf Proof.}  To convert from the Farey shift expansion~\numb{eq204} to
the additive continued fraction expansion~\numb{eq202}, we use
$\hc{D}=\hc{F}\hc{R}\hc{F}$ and $\hc{F}^2=\hc{I}$.

A finite automaton which converts the additive continued fraction
expansion~\numb{eq202} to the Farey tree expansion~\numb{eq204} must
keep two states to keep track of the sign $\pm 1$ of
$\det(\hc{S}_1\cdots\hc{S}_m)$ of the symbols $\hc{S}_i=\hc{F}$ or
$\hc{R}$ examined so far.  The initial state is $+1$.  $\qed$

The discussion in section 3.4 yields:

\begin{thm}~\label{MGCFtoCut}
For irrational $\theta$ with $- \frac{1}{2} < \theta < \frac{1}{2}$
there exists a finite automaton to
convert the Minkowski geodesic continued fraction expansion of $\theta$ to
that of the vertical cutting sequence expansion of $\langle \infty,
\theta \rangle$ and vice-versa.
\end{thm}

{\bf Proof.}  This follows from Theorem~\ref{MGCFtoCS}.  For each
direction, the finite automaton constructed needs
two states, to keep track of whether there a factor of $\hc{K}$ present;
i.e., to keep track of the sign of $\det(\hc{S}_1\cdots\hc{S}_n)$.  The
initial state is $+1$.
$\qed$

We illustrate two of these automata in Figure~3.1.
The
edge symbol $S: W$ specifies that this edge is taken if the current input
symbol is $S$,
and $W$ denotes a symbol sequence to be output.
The initial states are the vertices labeled $+1$.
\iffigures
\begin{figure}
%Figure 3.1
\setlength{\unitlength}{2in}
\centerline{ % Center this figure
\begin{picture}(2,.7)(-.5,-.3)
\put(0,0){\circle{.2}} % circles for automaton states
\put(1,0){\circle{.2}}
\put(0,0){\vclap{\clap{$+1$}}} % Labels of states
\put(1,0){\vclap{\clap{$-1$}}}
% Radii of circles are chosen so that arrowheads will be approximately
% aligned with the ends of the arcs
\put(.5,-1.6){\arc(.4,1.628){27.610}} % Arcs for transformations
\put(.5,1.6){\arc(-.4,-1.628){27.610}}
\put(1.15,0){\arc(-.112,-.1){276.379}}
\put(-.15,0){\arc(.112,.1){276.379}}
\put(.9,.028){\vector(4,-1){0}} % Arrowheads on arcs
\put(.1,-.028){\vector(-4,1){0}}
\put(1.038,-.1){\vector(-1,1){0}}
\put(-.038,-.1){\vector(1,1){0}}
\put(.5,.11){\clap{$\hc{F}:\emptyset$}} % Labels on arcs
\put(.5,-.16){\clap{$\hc{F}:\emptyset$}}
\put(1.15,.2){\clap{$\hc{R}:\hc{D}$}}
\put(-.15,.2){\clap{$\hc{R}:\hc{R}$}}
\end{picture}
}
\centerline{(a) Additive continued fraction to Farey tree}
\centerline{ % Center this figure
\begin{picture}(2,.8)(-.5,-.3)
\put(0,0){\circle{.2}} % circles for automaton states
\put(1,0){\circle{.2}}
\put(0,0){\vclap{\clap{$+1$}}} % Labels of states
\put(1,0){\vclap{\clap{$-1$}}}
\put(.5,-1.6){\arc(.4,1.628){27.610}} % Arcs for transformations
\put(.5,1.6){\arc(-.4,-1.628){27.610}}
\put(.5,-.4){\arc(.42,.483){81.980}}
\put(.5,.4){\arc(-.42,-.483){81.980}}
\put(1.15,0){\arc(-.112,-.1){276.379}}
\put(-.15,0){\arc(.112,.1){276.379}}
\put(.9,.028){\vector(4,-1){0}} % Arrowheads on arcs
\put(.1,-.028){\vector(-4,1){0}}
\put(.92,.083){\vector(4,-3){0}}
\put(.08,-.083){\vector(-4,3){0}}
\put(1.038,-.1){\vector(-1,1){0}}
\put(-.038,-.1){\vector(1,1){0}}
\put(.5,.27){\clap{$\hc{J}:\bar{\hc{J}}$}}
\put(.5,.11){\clap{$\hc{L}:\bar{\hc{R}}$}}
\put(.5,-.045){\clap{$\hc{L}:\bar{\hc{L}}$}}
\put(.5,-.205){\clap{$\hc{J}:\bar{\hc{J}}$}}
\put(1.15,.2){\clap{$\hc{R}:\bar{\hc{L}}$}}
\put(-.15,.2){\clap{$\hc{R}:\bar{\hc{R}}$}}
\end{picture}
}
\centerline{(b) MGCF to cutting sequence}
\caption{Finite Automata--Transducers (initial state for both automata
is $+1$).}
\end{figure}
\else %\iffigures
$$
\begin{array}{c}
\hline
~~~~~~~~~~~\mbox{Figure~3.1 about here}~~~~~~~~~ \\
\hline
\end{array}
$$
\fi

\section{Vertical Cutting Sequences to Additive Continued Fractions}
\hsp
The ordinary continued fraction expansion of an irrational real number
$ -  \frac{1}{2}  <  \theta  <  \frac{1}{2}$ is easily determined from
either its cutting
sequence expansion
or its Minkowski geodesic continued fraction expansion.
We begin with the latter case.

\begin{thm}\label{MGCFtoOCF}
Given $\theta$ with $-\frac{1}{2}   <  \theta  <  \frac{1}{2}$, with
ordinary
continued fraction expansion $\theta =  [a_0; a_1 , a_2 ,  a_3 , \cdots ]$,
let $\hc{S}_0 \hc{S}_1 \hc{S}_2 \cdots$ be the Minkowski geodesic
continued fraction expansion of $\theta$
in the alphabet $\{  \hc{L}  ,  \hc{R} , \hc{J}  ,  \hc{C}  \}$.
The MGCF expansion can be uniquely factored into segments $\hc{B}_0 \hc{B}_1
\cdots$
where $ \hc{B}_0 =  \hc{J}$ or $\hc{J}  \hc{L}$ and each $\hc{B}_i $ for $i
\geq 1$ is
$\hc{R}^k \hc{J}$ for some $k \geq 1$,
$\hc{R}^{k+1}\hc{J} \hc{L}$ for some $k \geq 1$,
or $\hc{R}^k \hc{C}$ for some $k \ge 1$.
Each segment encodes exactly one or two partial quotients of the OCF expansion.
For general segments, $\hc{R}^k \hc{J}$ encodes a partial quotient
$a_n =  k$,
while $\hc{R}^{k+1}\hc{J} \hc{L}$ or $\hc{R}^k \hc{C}$ encode two partial
quotients
$a_n =  k$,
$a_{n+1} =  1$.
For the segment $\hc{B}_0$, the symbol $\hc{J}$ encodes $a_0 =  0$,
while $\hc{J}  \hc{L}$ encodes the partial quotients $a_0 =  -1$, $a_1 =  1$.
The case $\hc{R}^k \hc{C}$ can only occur when $\theta$ is rational.
\end{thm}

{\bf Proof.}  We say that $p/q$ is a {\em best approximation} to
$\theta$ if $|q\theta-p|<|q'\theta-p'|$ for $0<q'<q$.  Using this
definition, it is a basic result about ordinary continued fractions that
the convergents are the complete set of best approximations; see Hardy
and Wright~\cite[Theorem 182]{HarWri60}.  We say that $p/q$ is a {\em
better approximation} to $\theta$ than $p'/q'$ if
$|q\theta-p|<|q'\theta-p'|$.  The best approximation property thus proves
that $p_n/q_n$ is a better approximation than any fraction with
denominator less than $q_{n+1}$.

We start with  $\hc{P}^{(0)} = \twobytwo 1 0 0 1 $, the identity matrix.
Thus $\hc{P}^{(0)} \hc{B}_t(\theta)$ is $\twobytwo 1 0 {-\theta} t $.
Neither of the Minkowski inequalities \numb{mink2} and
\numb{mink3} can hold with equality for any $t$, since
$|\theta|<\frac{1}{2}$.   We get $t_1^2
= 1 - {\left(\thfl - \theta\right)}^2$ as the value at which
\numb{mink1} holds with equality.  The first partial quotient matrix is
thus $\hc{J}$.  That makes $\hc{P}^{(1)}=\twobytwo 0 1 1 0 $.  Now
$\hc{P}^{(1)} \hc{B}_t(\theta) = \twobytwo {-\theta}
t 1 0 $.

We check the Minkowski inequalities.  Inequality \numb{mink1} is
satisfied, with equality only at $t=t_1$.  Inequality \numb{mink2} is an
equality if $1=(1 - \theta)^2 + t_2^2$, and \numb{mink3} is an equality
if $1=( - 1 - \theta)^2 + t_2^2$.  If $\theta=0$, neither one holds with
equality for any $t>0$, and $\hc{J}$ is thus the whole MGCF, encoding
$\theta=[0]$.  If $0<\theta<\frac{1}{2}$, then only \numb{mink2} can
hold with equality, so the next partial quotient matrix will be
$\hc{R}$.  In this case, we have
\[ \hc P^{(1)} = \twobytwo 0 1 1 0
               = \twobytwo {p_0}{q_0}{p_{-1}}{q_{-1}}, \]
and we let $\hc{B}_0=\hc{J}$.
If $-\frac{1}{2}<\theta<0$, then only \numb{mink3} can hold with equality,
so the next partial quotient matrix will be $\hc{L}$.
For $-\frac{1}{2}<\theta<0$, the continued fraction for $\theta$ begins
$[-1, 1, \ldots\,]$, and thus $p_1 = 0, q_1 = 1$.
This gives
\[ \hc{P}^{(2)} = \twobytwo 0 1 {-1} 1
               = \twobytwo {p_1}{q_1}{p_0}{q_0}, \]
and we let $\hc{B}_0=\hc{J},\hc{L}$.

In either case, after the first step, we have encoded all coefficients
through $a_{n-1}$, and our current matrix is
\BE
\hc{P}^{(j)}=\twobytwo
{p_{n-1}}{q_{n-1}}{p_{n-2}}{q_{n-2}}. \label{stdform}
\EE
The rest of the theorem is proved by induction, with the induction
hypothesis that after each segment $\hc{B}_k$, the matrix $\hc{P}^{(j)}$
is in the form \numb{stdform}, with the coefficients through $a_{n-1}$
encoded.
\smallbreak

First, assume that the current $\hc{P}^{(j)}$ is in the form
\numb{stdform}.
Then
\[ \hc{P}^{(j)} \hc{B}_t(\theta) =
   \twobytwo {\app{n-1}}{tq_{n-1}}{\app{n-2}}{tq_{n-2}}. \]
Since $q_{n-1}>q_{n-2}$, the requirement of an increasing denominator
makes $\hc{J}$ impossible for the next partial quotient.
The Minkowski inequality for $\hc{R}$ is
\begin{eqnarray} % {} needed before plus sign to get proper spacing
   {\bigl(\papp{n-1}+\papp{n-2}\bigr)}^2 - {\papp{n-2}}^2 \nonumber \\
      {}+ t^2 \bigl((q_{n-1}+q_{n-2})^2 - q_{n-2}^2\bigr) & \ge & 0,
\label{stdmink2}
\end{eqnarray}
and for $\hc{L}$, it is
\begin{eqnarray} % {} needed before plus sign to get proper spacing
   {\bigl(\papp{n-1}-\papp{n-2}\bigr)}^2 - {\papp{n-2}}^2 \nonumber \\
      {}+ t^2 \bigl((q_{n-1}-q_{n-2})^2 - q_{n-2}^2\bigr) & \ge & 0.
\label{stdmink3}
\end{eqnarray}
By the best approximation condition, $p_{n-1}/q_{n-1}$ is a better
approximation to $\theta$ than $(p_{n-1} - p_{n-2}) \big/ (q_{n-1} -
q_{n-2})$.  Thus~\numb{stdmink3} is satisfied for small $t$, and thus
for all $t<t_j$, so we cannot have $\hc{L}$ next.  If $\theta$ is
rational and $a_{n-1}$ is the last term, then~\numb{stdmink2} is
satisfied for all $t$; thus the MGCF terminates.  Otherwise, since every
even convergent is less than $\theta$, and every odd convergent is
greater, $\app{n-1}$ and $\app{n-2}$ have opposite signs, and the
absolute value of their sum is less than the absolute value of
$\app{n-2}$.  Thus~\numb{stdmink2} is not satisfied for sufficiently
small $t$, and $t_{j+1}$ is the value of $t$ at which it holds with
equality.  Thus the next matrix is $\hc{R}$; hence
\[ \hc{P}^{(j+1)} = \hc{R}\hc{P}^{(j)} = \twobytwo
{p_{n-1}}{q_{n-1}}{p_{n-2}+p_{n-1}}{q_{n-2}+q_{n-1}}. \]
\smallbreak

The remainder of the proof considers the case where $\hc{P}^{(j)}$ is of
the form
\BE
\hc{P}^{(j)}=
\twobytwo {p_{n-1}}{q_{n-1}}{p_{n-2}+mp_{n-1}}{q_{n-2}+mq_{n-1}},
\label{intform}
\EE
which can be reached after $m$ applications of $\hc{R}$ to the
form~\numb{stdform}.  Now we have
\[ \hc{P}^{(j)}\hc{B}_t(\theta) =
   \twobytwo {\app{n-1}}{tq_{n-1}}
	     {\papp{n-2}+m\papp{n-1}}{t(q_{n-2}+mq_{n-1})}.
\]
Since $q_{n-2}+mq_{n-1}>q_{n-1}$, the next symbol $A^{(j+1)}$ cannot be
$\hc{L}$.  The other two Minkowski
inequalities, for $\hc{J}$ and $\hc{R}$, are
\begin{eqnarray} % {} needed before plus sign to get proper spacing
   {\bigl(\papp{n-2}+m\papp{n-1}\bigr)}^2 - {\papp{n-1}}^2 \nonumber \\
      {}+ t^2 \bigl((q_{n-2}+mq_{n-1})^2 - q_{n-1}^2) & \ge & 0,
         \label{dowehit}\\
   {\bigl((m+1)\papp{n-1}+\papp{n-2}\bigr)}^2 \nonumber \\
      {}- \bigl(m\papp{n-1}+\papp{n-2}\bigr)^2 \nonumber \\
      {}+ t^2 \bigl(((m+1)q_{n-1} + q_{n-2})^2-(mq_{n-1} + q_{n-2})^2\bigr)
      & \ge & 0. \label{dowemiss}
\end{eqnarray}
Thus the next symbol is $\hc{A}^{(j+1)}=\hc{J}$ if~\numb{dowehit} holds
with equality for a larger $t<t_j$ than~\numb{dowemiss},
$\hc{A}^{(j+1)}=\hc{R}$ if~\numb{dowemiss} holds with equality for a
larger $t<t_j$, and $\hc{A}^{(j+1)}=\hc{C}$ if both hold with equality
at the same $t$.  Recall that we have, for $0\le m \le a_n$, and $n$ even,
\BE   \frac{p_{n-1}+p_n}{q_{n-1}+q_n}
    = \frac{(a_n+1)p_{n-1}+p_{n-2}}{(a_n+1)q_{n-1}+q_{n-2}}
    > \theta \ge \frac{mp_{n-1} + p_{n-2}}{mq_{n-1} + q_{n-2}}, \EE
with the reverse inequalities holding for $n$ odd.  The last inequality
holds with equality only if $\theta = p_n/q_n$ and $m=a_n$.
If $m<a_n$, then $mq_{n-1} + q_{n-2} < q_n$, so we must have
\[
\papp{n-1}^2 < \bigl(\papp{n-2}+m\papp{n-1}\bigr)^2,
\]
and thus the left-hand side of~\numb{dowehit} is positive for small $t$,
so $\hc{A}^{(j+1)}$ cannot be $\hc{J}$ or $\hc{C}$.  Thus, if $m<a_n$, the
only possible case is $\hc{R}$, and the matrix $\hc{P}^{(j+1)}$ is still
of the form~\numb{intform}, with $m$ replaced by $m+1$.

We next consider the case when $m=a_n$ in~\numb{intform}.  In this case, the
next symbol $\hc{A}^{(j+1)}$ may be any of $\hc{R}$, $\hc{J}$, and
$\hc{C}$.  First suppose $\hc{A}^{(j+1)}=\hc{J}$.  If so, we
have finished our segment $\hc{B}_k=\hc{R}^{a_n}\hc{J}$, and
$\hc{P}^{(j+1)}$ is now in the correct form~\numb{stdform}, since it is
\[ \twobytwo {p_n}{q_n}{p_{n-1}}{q_{n-1}} . \]

Next suppose that $m=a_n$ and $\hc{A}^{(j+1)}=\hc{R}$.  For this to
happen, we
need \numb{dowemiss} to be true only for $t\ge t_c$.
Since $m=a_n$, we have $mq_{n-1} + q_{n-2} = q_n$, and similarly for
$p_n$.  For \numb{dowemiss} to be true only for sufficiently large
$t$, $(p_{n-1} + p_n)/(q_{n-1} + q_n)$ must be a better approximation to
$\theta$ than $p_n/q_n$.
This cannot happen if $\theta = p_n/q_n$, hence
$a_{n+1}$ is defined.  By the best approximation property, we must have
$q_{n-1}+q_n \ge q_{n+1}$, and it follows that
$a_{n+1}=1$, $p_{n+1}=p_n+p_{n-1}$, and  $q_{n+1}=q_n+q_{n-1}$, so that
in this case
\[
\hc{P}^{(j+1)} = \hc{R}\hc{P}^{(j)} =
\twobytwo {p_{n-1}}{q_{n-1}}{p_{n+1}}{q_{n+1}} .
\]
This matrix is in the form~\numb{intform} but with $m=a_n+1$.  Now
$\hc{A}^{(j+2)}$ cannot be $\hc{R}$ or $\hc{C}$, because the  left-hand
side of~\numb{dowemiss} is positive for small $t$, so
$\hc{A}^{(j+2)}$ must be $\hc{J}$, which gives
\BE
\hc{P}^{(j+2)}= \hc{J}\hc{R}\hc{P}^{(j)} =
\twobytwo {p_{n+1}}{q_{n+1}}{p_{n-1}}{q_{n-1}} . \label{missform} \EE
This gives
\[ \hc{P}^{(j+2)}\hc{B}_t(\theta) =
\twobytwo {\app{n+1}}{tq_{n+1}}{\app{n-1}}{tq_{n-1}} , \]
The following symbol $\hc{A}^{(j+3)}$ cannot be $\hc{J}$, because
$q_{n+1}>q_{n-1}$ and we must have
an increasing denominator.  The Minkowski inequality for $\hc{R}$ is
\begin{eqnarray} % {} needed before plus sign to get proper spacing
   {\bigl(\papp{n+1}+\papp{n-1}\bigr)}^2 - {\papp{n-1}}^2 \nonumber \\
      {}+ t^2 \bigl((q_{n+1}+q_{n-1})^2 - q_{n-1}^2\bigr) & \ge & 0.
\end{eqnarray}
The constant term here is positive because $\app{n+1}$ and $\app{n-1}$
have the same sign, so this inequality holds for all $t$.  The Minkowski
inequality for $\hc{L}$ is
\begin{eqnarray} % {} needed before plus sign to get proper spacing
   {\bigl(\papp{n+1}-\papp{n-1}\bigr)}^2 - {\papp{n-1}}^2 \nonumber \\
      {}+ t^2 \bigl((q_{n+1}-q_{n-1})^2 - q_{n-1}^2\bigr) & \ge & 0,
\end{eqnarray}
which is not satisfied for small enough $t$ because
$\papp{n+1}-\papp{n-1}=\app{n}$, and thus its constant term is negative.
Thus $\hc{A}^{(j+3)}=\hc{L}$, so that
\BE
\hc{P}^{(j+3)} = \hc{L}\hc{J}\hc{R}\hc{P}^{(j)} =
\twobytwo {p_{n+1}}{q_{n+1}}{p_n}{q_n}.
\EE
This matrix has the form \numb{stdform}, and we have encoded the two
coefficients, $a_n$ and $a_{n+1}=1$, with a segment
$\hc{B}_k=\hc{R}^{a_{n}+1} \hc{J} \hc{L}$.

Finally, suppose that $m=a_n$ and
$\hc{A}^{(j+1)}=\hc{C}$.  Again, \numb{dowemiss} must fail to hold for
sufficiently small $t$, and it follows that $a_{n+1}=1$.  Thus,
since $q_{n+1} = q_n +
q_{n-1}$ and $p_{n+1} = p_n + p_{n-1}$, we have
\BE
 \hc{P}^{(j+1)} = \hc{C} \hc{P}^{(j)} =
  \twobytwo 1 1 0 1 \twobytwo {p_{n-1}}{q_{n-1}}{p_n}{q_n}
   = \twobytwo {p_{n+1}}{q_{n+1}}{p_n}{q_n}.
\EE
Here $\hc{P}^{(j+1)}$ is in the form~\numb{stdform}, and we have
again encoded the two coefficients, $a_n$ and $a_{n+1} = 1$ as
$\hc{R}^{a_n}\hc{C}$.  This completes the induction step.

Lemma~\ref{vertcorner} shows that the last case $\hc{R}^k\hc{C}$
can occur only for rational $\theta$. $\qed$

There is an analogous conversion method from the cutting sequence
expansion to the additive ordinary continued fraction expansion, as
follows.

\begin{thm}\label{CStoOCF}
Given $\theta$ with ordinary continued function expansion $\theta = [
a_0; a_1 , a_2 , a_3 , \cdots ]$, let $\hc{S}_0^\ast \hc{S}_1^\ast
\hc{S}_2^\ast \cdots$ be the cutting sequence expansion for the geodesic
$\langle \In , \theta \rangle = \{ \theta + it : t > 0 \}$ in the
alphabet $\{ \bar{\hc{L}} , \bar{\hc{R}} , \bar{\hc{J}} , \bar{\hc{C}}_1
, \bar{\hc{C}}_2 \}$.  It can be uniquely factored into segments
$\bar{\hc{B}}_0 \bar{\hc{B}}_1 \bar{\hc{B}}_2 \cdots$ where
$\bar{\hc{B}}_0 = \bar{\hc{J}} $ or $\bar{\hc{J}} \bar{\hc{R}} $,
encoding $a_0=0$ or $a_0=-1$, $a_1=1$, respectively, and each succeeding
segment is $\bar{\hc{R}}^k \bar{\hc{J}} $ or $\bar{\hc{L}}^k \bar{\hc{J}}
$ encoding $a_n$ for $n$ even and odd, respectively, or is
$\bar{\hc{R}}^{k+1} \bar{\hc{J}} ~\bar{\hc{R}}$, or $\bar{\hc{R}}^k
\bar{\hc{C}}_1 $ encoding $a_n = k , a_n+1 = 1$ for $n$ even, or is
$\bar{\hc{L}}^{k+1} \bar{\hc{J}} ~ \bar{\hc{L}}$ or
$\bar{\hc{L}}^k\bar{\hc{C}}_2 $ encoding $a_n = k$, $ a_{n+1} = 1$ for
$n$ odd.  The symbols $\bar{\hc{C}}_1 , \bar{ \hc{C}}_2 $ can only occur
in expansions of rational $\theta$.

\end{thm}

{\bf Proof.}  This follows from Theorem \ref{MGCFtoOCF} by noticing that
the parity of the initial word $\hc{W}_i=\hc{B}_0\ldots \hc{B}_i$
changes after each segment $\hc{B}_i$ encoding one term, and does not
change after any segment $\hc{B}_i$ encoding two terms, while the
segment $\hc{B}_0$ has odd parity if it encodes no terms, and even
parity if it encodes $a_1$.  $\qed$

An immediate consequence of Theorem \ref{MGCFtoOCF} and \ref{CStoOCF} is
the following result.

\begin{thm}\label{autMGCFtoOCF}
There exists a finite automaton which converts the Minkowski
geodesic continued fraction expansion of each irrational $\theta$
with $- \frac{1}{2} < \theta < \frac{1}{2}$ to the additive continued
fraction expansion of $\theta$.
There exists a finite
automaton that converts the vertical cutting sequence expansion of
$\langle \In , \theta \rangle$ to the additive
ordinary continued fraction expansion of $\theta$.
\end{thm}

{\bf Proof.}
The segment-partition of the MGCF expansion given in
Theorem~\ref{MGCFtoOCF} permits the Farey tree
expansion (\ref{eq204}) to be computed by a finite automaton, because
the necessary information to decide on the symbol $\hc{R}$ versus
$\hc{L}$ depends only on the determinant $ \pm 1$ of the product of the
MGCF matrices scanned plus the values of the last four MGCF symbols in
the expansion.  Next, Theorem~\ref{CStoOCF} guarantees that a finite
automaton exists to convert the cutting sequence expansion as well.
Finally, Theorem~\ref{ACFtoFarey} applies to give the additive continued
fraction from the Farey tree expansion.
$\qed$

Theorem \ref{MGCFtoOCF} implies that the Minkowski geodesic continued
fraction expansion of $\theta$ can represented in an abbreviated form
$$
\theta  :  = [ \tilde{a}_0,\tilde{a}_1 ,\tilde{a}_2 , \tilde{a}_3 , . . . ]
$$
similar to its ordinary continued fraction expansion
$$
\theta =  [a_0; a_1 , a_2 ,  a_3 , . . . ] ~,
$$
with the change that each symbol~1 in the OCF expansion is to be
replaced by one of three possible symbols $1_h$, $1_m$ and $1_c$.  Here
$1_h $ means that the continued fraction partial quotient $a_n = 1$
begins a new segment $\hc{R}\hc{J}$ or $\hc{R}\hc{R}\hc{J}\hc{L}$ in
the MGCF expansion (so that the previous convergent $p_{n-1}/q_{n-1}$
was ``hit'' at the end of a segment), $1_m$ means that it combines with
the previous partial quotient in a block $\hc{R}^{k+1}\hc{J}\hc{L}$ (so that
$p_{n-1}/q_{n-1}$ was ``missed''), and $1_c$ means that it combines with
the previous partial quotient in a block with a $\hc{C}$-symbol
(a ``corner'').  However it does not specify which
sequence of symbols actually occur as legal expansions.  We study this
next.

\section{Additive Continued Fractions to Vertical Cutting Sequences}
\hsp
In this section, we show how to construct the Minkowski geodesic continued
fraction
expansion of $\theta$ given the additive continued fraction expansion of
$\theta$.
Let $\theta$ have the ordinary continued fraction expansion
$$
\theta =  [a_0; a_1 , a_2 ,  a_3 , \cdots ],
$$
in which $a_0=0$ or $-1$.
In view of the results of section 4, it suffices to determine for each $a_n
=  1$ whether or not
it has label $1_h  , 1_m$ or $1_c$ in the partition of the MGCF
expansion given in Theorem \ref{MGCFtoOCF}.

\begin{thm}\label{1h1mthm}
Given $\theta =  [a_0; a_1 , a_2 ,  a_3 , \cdots ]$ with
$-\frac{1}{2}  <  \theta  <  \frac{1}{2}$,
suppose $a_{n+1} =  1$.
Set
$\alpha_n   =  [ 0, a_n , a_{n-1}  \dd  a_1 ] = \df{q_{n-1}}{q_n} $,
and
$\beta_n   =   [a_{n+1} , a_{n+2} , \cdots ] $,
so that $\alpha_n \in [0,1]$ and $\beta_n \in [1,2]$.
In terms of the linear fractional transformation
\beql{eq500}
\hc{N} (z) := \left[ \begin{array}{cc}
1 & 2 \\
2 & 1
\end{array}
\right] (z) = \df{z+2}{2z+1}~,
\eeq
the Minkowski geodesic continued fraction expansion of $\theta$ has
\beql{1h1mcond}
\tilde{a}_{n+1} =  \left\{
\begin{array}{ll}
1_h ~~~ & \mbox{if ~$\beta_n > \hc{N}( \alpha_n )$} ~, \\
1_c ~~~ & \mbox{if ~$\beta_n = \hc{N}( \alpha_n ) $}~, \\
1_m ~~~ & \mbox{if ~$\beta_n < \hc{N}( \alpha_n )$}~.
\end{array}
\right.
\eeq
\end{thm}

{\bf Remark.} The matrix $\hc{N}$ acting as a linear fractional
transformation maps $[0,1]$ to $[1, 2]$ while reversing orientation.
Its inverse $\hc{N}^{-1}=\frac{1}{3}\twobytwo{-1}{2}{2}{-1}$ is not
integral.
The theorem could also be formulated in terms of the linear fractional
transformation
$-\hc{N}$ which maps $[0,1]$ to $[-2,-1]$ and is an
involution, i.e. $-  \hc{N}  ( -  \hc{N}  (z))  \equiv  z$.
In particular,
since $\hc{N} $ sends $\RR^+$ into itself,
\beql{eq404}
\beta  <  \hc{N} ( \alpha ) \iff \alpha  >  - \hc{N} ( - \beta ) ~.
\eeq
The symbol $1_c$ can occur only when $\theta$ is rational, because
$\hc{N} ( \alpha_n )$ is always rational,
while $\beta_n$ is rational if and only if the expansion of $\theta$
terminates.

{\bf Proof of Theorem \ref{1h1mthm}.}
Note that if $a_1 = 1$, it always becomes $1_m$, and we have $\alpha_n=0$
and thus $\hc{N}(\alpha_n) = 2 > \beta_n$ as required. The
following discussion assumes that $a_{n+1} = 1$, with $n \ge 1$, so that
we are not in the segment $\hc{B}_0$.

By the argument in the proof of Theorem~\ref{MGCFtoOCF}, whether
$a_{n+1}$ becomes $1_h$, $1_c$, or $1_m$ depends on which of
\numb{dowehit} and \numb{dowemiss} holds with equality for a larger
value of $t$; we have $1_m$ if it is \numb{dowemiss}, and $1_c$ if both
hold with equality for the same $t$.  Since $a_{n+1} = 1$, we have
\begin{eqnarray*}
   p_{n+1}  & = & p_n + p_{n-1}. \\
   q_{n+1}  & = & q_n + q_{n-1}. \\
   \app{n+1} & = & \papp{n} + \papp{n-1},
\end{eqnarray*}
The condition for $1_m$ is thus
\BE \frac{\papp{n}^2 - \papp{n+1}^2}{q_{n+1}^2-q_n^2}
   >\frac{\papp{n-1}^2 - \papp{n}^2}{q_n^2-q_{n-1}^2}, \EE
or, equivalently,
\BE \frac{\papp{n}^2 - \papp{n+1}^2}{\papp{n-1}^2 - \papp{n}^2}
   >\frac{q_n^2 - q_{n+1}^2}{q_{n-1}^2 - q_n^2}, \label{fracmiss} \EE
and the condition for $1_c$ is equality.  Putting everything in terms of
$p_{n-1}$, $p_n$, $q_{n-1}$, and $q_n$ gives
\BE \frac{\papp{n-1}^2 + 2\papp{n-1}\papp{n}}{\papp{n-1}^2 - \papp{n}^2}
   >\frac{q_{n-1}^2 + 2q_{n-1}q_n}{q_{n-1}^2 - q_n^2}. \label{intmiss} \EE
Now let
\begin{eqnarray}
   \alpha_n & = & \frac{q_{n-1}}{q_n}, \\
   \beta_n & = & -\frac{\app{n-1}}{\app{n}}.
\end{eqnarray}
These quantities have ordinary continued-fraction expansions
$\alpha_n   =  [ 0, a_n , a_{n-1}  \dd  a_1 ]$,
and
$\beta_n   =   [a_{n+1} , a_{n+2} , \cdots ] $; see Venkov~\cite[section
2.4] {Ven70}.
Clearly $0 < \alpha_n < 1$ and $\beta_n > 1$, and, since $a_{n+1}=1$,
$\beta_n < 2$.
Now \numb{intmiss} becomes
\BE \frac{\beta_n^2 - 2\beta_n}{\beta_n^2 - 1}
   >\frac{\alpha_n^2 + 2\alpha_n}{\alpha_n^2 - 1}, \EE
or, by subtracting 1 from each side,
\BE \frac{-2\beta_n + 1}{\beta_n^2 - 1}
   >\frac{2\alpha_n + 1}{\alpha_n^2 - 1}. \label{quadmiss} \EE
Let $C = (2\alpha_n + 1)/(\alpha_n^2 - 1)$.  Then
\numb{quadmiss} holds if and only
if $-\beta_n$ is between the roots of
\BE Cx^2 - 2x - 1 - C = 0. \label{quadeqn} \EE
The sum of the roots of \numb{quadeqn} is $2/C$, and the larger root is
$\alpha_n$, so the other root is
\BE
x = \frac{2\alpha_n^2-2}{2\alpha_n+1}-\alpha_n =
\frac{-\alpha_n-2}{2\alpha_n+1} = -\hc{N}(\alpha_n).
\EE
This gives our desired condition; we have $\tilde{a}_{n+1}=1_m$ if
$\beta_n<\hc{N}(\alpha_n)$.  Since $-\hc{N}$ is an involution and is
increasing on  $[-2,-1]$ and on $[0,1]$, we can also write this
condition as $\alpha_n>-\hc{N}(-\beta_n)$.
$\qed$

Theorem \ref{1h1mthm} suffices to classify all $\theta$ containing the symbol
$1_c$, i.e. all vertical geodesics that hit a corner of a translate of a
fundamental domain.

\begin{cor}\label{autfinds1c}
Let $-\frac{1}{2}  <  \theta  <  \frac{1}{2}$ and suppose that the
MGCF expansion of $\theta$ contains a symbol $1_c$.
Then $\theta$ is rational and
has ordinary continued fraction expansion of the form
\beql{eq405}
\theta =
\cases{
[ 0, a_1, a_2,  \dd  a_n , 1_c , b_1  \dd  b_m ],
  & if $0\le\theta<\frac{1}{2}$.\cr
[ -1, 1, a_1-1, a_2,  \dd  a_n , 1_c , b_1  \dd  b_m ],
  & if $-\frac{1}{2}<\theta<0$.\cr
}
\eeq
in which the additive continued fraction expansion of $\alpha_n^*  :  =
[1, b_1  \dd  b_m ]$ is computable from the additive continued fraction
expansion of $\alpha_n =  [0, a_n  \dd  a_1 ]$ by a finite automaton,
and vice-versa.
Furthermore there is an absolute constant $c_0 $ such that
\beql{eq406}
\mbox{ ACF-length ~$ ( \alpha_n^* )  \leq  c_0$  ~ ( ACF-length~ $(
\alpha_n )) $}~.
\eeq
\end{cor}

{\bf Proof.}
Raney \cite{Ran73} proves that given any fixed linear fractional
transformation $\tilde{\hc{M}} (z) : = \frac{az+b}{cz+d}$ with integer
coefficients $a,b,c,d$ and $ad-bc \neq 0$, the ACF expansion of
$\tilde{\hc{M}} ( \theta )$ can be computed from the ACF expansion of
$\theta$ using a finite automaton.  Furthermore the conversion process
inflates the ACF-length by at most a multiplicative constant (depending on
$\tilde{\hc{M}}$).  Apply this to $\beta_n = \hc{N} ( \alpha_n )$.  The
specific automaton for $\hc{N}(z)=\frac{z+2}{2z+1}$ is given
in Raney~\cite[pp.~274--275]{Ran73}, and the constant in this case is 3.
$\qed$

These results give a bound on the computational complexity of
computing the MGCF expansion from the
additive continued fraction expansion.
\begin{thm}\label{turing}
The Minkowski geodesic continued fraction expansion of $\theta$ can be
computed from the additive continued fraction expansion in quadratic
time using linear space.  That is, there are absolute constants $c_1,
c_2, c_3$ such that for any $\theta$ with
$-\frac{1}{2}<\theta<\frac{1}{2}$, the first $\ell$ symbols of the MGCF
expansion of $\theta$ can be computed using the first $c_1 \ell$ symbols
of the additive continued fraction expansion of $\theta$ using at most
$c_2 \ell^2$ time steps and $c_3 \ell$ space locations.
\end{thm}

{\bf Proof.}  Theorem~\ref{MGCFtoOCF} shows that the main problem is to resolve
whether a given symbol $a_{n+1}=1$ which appears as the $\ell$th symbol
in the ordinary continued fraction expansion of $\theta$ is to be $1_m$,
$1_h$, or $1_c$ in the MGCF expansion.  For this, we use
Theorem~\ref{1h1mthm}.  At worst, we must look all the
way back to the beginning of the MGCF expansion.
Corollary~\ref{autfinds1c} shows that the ACF expansion for
$\hc{N}(\alpha_n)$ is of length at most $c_0\ell$, and thus comparing it
to $\beta_n$ requires looking at no more than $c_0\ell+1$
symbols, and this uses $O(\ell)$ time and $O(\ell)$ space.  Since there may
be $O(\ell)$ different 1's to be resolved, the total time required is
$O(\ell^2)$. Note that we test the inequalities 
(\ref{1h1mcond})
by comparing symbol sequences for $\hc{N}(\alpha_n)$ with the initial
part of that for $\beta_n$.
This algorithm can easily be implemented on a random access machine
with the given time and space bounds.
It can also be implemented
on a one-tape Turing machine with the same space bound and
a time bound polynomial in $\ell$.  We omit details. For these 
two standard computational machine models, 
see Aho, Hopcroft and Ullman~\cite{AhoHopUll74}.
$\qed$

There are examples which do require $\Omega(\ell^2)$ time
steps; for example,
\[
\theta=(\sqrt{3}-1)/2=[0,\overline{2,1_h}],
\]
has this property.  For this $\theta$, it is necessary to backtrack all
the way to the first symbol to determine that each 1 is $1_h$,
because the sequence $[0,(2, 1_h )^j , 2, 1_c , (2, 1_m )^j ,4]$
is a $1_c$-sequence for each $j \geq 1$.

In the next section, we will show that the MGCF expansion cannot be
computed from the additive continued fraction expansion using a finite
automaton, or even using a pushdown automaton with one stack, as defined in
Hopcroft and Ullman~\cite{HopUll79}.

\section{Vertical Geodesics: Forbidden Blocks}
\hsp
In this section, we characterize the allowable cutting sequences of
$\Pi_\sF^0$.

\begin{defn}
A finite word $\hc{W}$ in the symbol set
$\{\bar{\hc{L}},\bar{\hc{R}},\bar{\hc{J}}\}$ is a {\em forbidden block}
of $\Pi_\sF^0$ if it occurs in no cutting sequence of $\Pi_\sF^0$;
otherwise, it is an {\em admissible block.}  It is an {\em excluded
initial block} of $\Pi_\sF^0$ if does not occur as an initial segment of
any cutting sequence of $\Pi_\sF^0$; otherwise, it is an {\em
included initial block}.
\end{defn}

It is easy to see that excluded initial blocks alone can be used to
characterize
{\em any} subset of the one-sided shift on
$\{\bar{\hc{L}},\bar{\hc{R}},\bar{\hc{J}}\}$.  We show that
$\Pi_\sF^0$ has the stronger property that it is
determined by its set of forbidden blocks,
as follows.

\begin{thm}\label{initblocks}
A finite word $\hc{W}$ in the symbols $\bar{\hc{L}},\bar{\hc{R}},\bar{\hc{J}}$
is an excluded initial block in $\Pi_\sF^0$ if and only if
at least one of the blocks $\bar{\hc{L}}\hc{W}$ and
$\bar{\hc{R}}\hc{W}$ is a forbidden block of $\Pi_\sF^0$.
\end{thm}

{\bf Proof.}  We prove the contrapositive.  First, we show
that for every included initial block $\hc{W}$, both
$\bar{\hc{R}}^n\hc{W}$ and $\bar{\hc{L}}^n\hc{W}$ are admissible blocks
of $\Pi_\sF^0$ for all $n\ge 1$.  Geometrically, this encodes the fact
that the geodesic $\langle\infty,\theta\rangle$ is a limit of geodesics
$\langle\theta',\theta\rangle$ as $|\theta'|\to\infty$, where
$\theta'\to-\infty$ is associated to $\bar{\hc{R}}^n\bar{\hc{W}}$ and
$\theta\to+\infty$ is associated to $\bar{\hc{L}}^n\bar{\hc{W}}$.  The
block $\hc{W}$ corresponds to a finite initial segment of the cutting
sequence of $\theta+it$ for some irrational $\theta$, say for $t_0\le
t\le\infty$ and it has first symbol $\bar{\hc{J}}$.  By
Lemma~\ref{vertcorner}, this geodesic cannot hit a corner of any
translate of $\sF$, hence there is some positive $\epsilon$ such that it is at
distance at least $\epsilon$ from any corner for $t_0\le t\le\theta$.
Pick $p/q$ a large negative rational number, and observe that the
geodesic $\langle\theta',\theta\rangle$ contains the word
$\bar{\hc{R}}^n\hc{W}$ in its one-sided infinite cutting sequence, in which
it enters the domain $\sF$ just after the $\bar{\hc{R}}^n$.  Here we
require that the radius $r=\theta-p/q$ is at least $n+2$ so that it
produces the sequence
$\bar{\hc{R}}^n$, and satisfies $r-\sqrt{r^2-1}<\epsilon$ so that the geodesic
is within $\epsilon$ of $\langle\infty,\theta\rangle$ over the range
$t_0\le t<1$ and thus produces the sequence $\hc{W}$.  Now take a matrix
$\hc{M}=\twobytwo{q'}{p'}{q}{-p}\in SL(2,\ZZ)$, so that
$\hc{M}(\gamma)=\langle\infty,\theta'\rangle$ with
$\theta'=\hc{M}(\theta)$, and note that $\theta'$ is necessarily
irrational.  The
$PSL(2,\ZZ)$-action doesn't affect cutting sequences, so the cutting
sequence of $\hc{M}(\gamma)$ still contains the word
$\bar{\hc{R}}^n\hc{W}$.  There is a unique choice of $q'$ and $p'$ such
that $-\frac{1}{2}<\theta'<\frac{1}{2}$.  This vertical geodesic has the
same cutting sequence as $\gamma$.  Thus $\bar{\hc{R}}^n\hc{W}$ and
likewise $\bar{\hc{L}}^n\hc{W}$ are admissible blocks of $\Pi_\sF^0$ for
all $n\ge 1$.

Second, we prove that if $\hc{W}=\hc{S}_1\cdots\hc{S}_n$ is a word in
$\bar{\hc{L}},\bar{\hc{R}},\bar{\hc{J}}$ for which $\bar{\hc{R}}\hc{W}$
and $\bar{\hc{L}}\hc{W}$ are both admissible blocks of $\Pi_\sF^0$,
there is a vertical geodesic whose cutting sequence has $\hc{W}$ as an
initial segment.  To show this, note that the first symbol in $\hc{W}$
is necessarily $\bar{\hc{J}}$, since $\bar{\hc{L}}\bar{\hc{R}}$ and
$\bar{\hc{R}}\bar{\hc{L}}$ are forbidden blocks.  Let $\gamma_1$ and
$\gamma_2$ be irrational vertical geodesics whose cutting sequences
contain the words $\bar{\hc{R}}\hc{W}$ and $\bar{\hc{L}}\hc{W}$,
respectively. Translate each of them by the appropriate elements of $PSL
(2,\ZZ)$ to geodesics $\gamma'_1=\langle p_1/q_1,\theta'_1\rangle$ and
$\gamma'_2=\langle p_2/q_2,\theta'_2\rangle$ in such a way that the
translated geodesics enter the fundamental domain $\sF$ at the symbol
immediately preceding $\hc{W}$ in their cutting sequences.  Since the
cutting sequence of $\gamma'_1$ has the letter $\bar{\hc{R}}$ as it
enters $\sF$, it is oriented in the direction of increasing real part
and has $p_1/q_1<-1/2$, while $\gamma'_2$ is oriented in the direction
of decreasing real part and has $p_2/q_2>1/2$.  The first symbol in
$\hc{W}$ is $\bar{\hc{J}}$, hence each geodesic exits $\sF$ through its
bottom edge at a point with real part between $-1/2$ and $1/2$.  If
$\gamma'_1$ and $\gamma'_2$ do not intersect, that is if
$\theta'_1<\theta'_2$, then pick any irrational $\theta$ between
$\theta'_1$ and $\theta'_2$, say $\theta=\theta'_1$, and let
$\gamma=\langle\infty,\theta\rangle$.  If $\gamma'_1$ and $\gamma'_2$ do
intersect, let $\theta$ be the real part of their intersection and let
$\gamma=\langle\infty,\theta\rangle$.  In either case, $\gamma$ passes
through $\sF$ since
$-\frac{1}{2}<\theta'_1\le\theta\le\theta'_2<\frac{1}{2}$.  Thus
the intersection must have real part $\theta$ with
$-\frac{1}{2}<\theta<\frac{1}{2}$.  In either case, we show that
$\gamma$ has initial word $\hc{W}$ in its cutting sequence.  After
leaving $\sF$,
both $\gamma'_1$ and $\gamma'_2$ pass through the same sequence of
translated fundamental domains $\sF_1,\ldots,\sF_n$, entering along the
same edge of each, corresponding to the symbols in $\hc{W}$.  The
vertical geodesic $\gamma$ must then enter $\sF_j$ on the same edge at a
point $z_j$ between the points $z_{1,j}$ and $z_{2,j}$ where $\gamma'_1$
and $\gamma'_2$ hit it, because the domains $\sF_j$ are hyperbolically
convex; see Figure 6.1.  Also, it cannot hit any other fundamental
domain between $z_{j-1}$ and $z_j$, because every point in that interval
on $\gamma$ is on a geodesic between points of $\gamma'_1$ and
$\gamma'_2$ which are in $\sF$.  Thus the initial cutting sequence of
$\gamma$ is $\hc{W}$.  Finally, if the finite endpoint $\theta$ of
$\gamma$ is rational, then since $\gamma$ hits no corners up to hitting
the edge $\hc{W}_n$, it is a Euclidean distance at least $\epsilon$ away
from every corner on the $\hc{W}$-edges.  For every irrational $\theta'$
with $|\theta'-\theta|<\epsilon$, the irrational vertical geodesic
$\gamma'=\{\theta'+it: t>0\}$ has the same initial segment $\hc{W}$.
$\qed$

\iffigures
\begin{figure}
% Figure 6.1
% Using domain bounded by (1/5,1/3), (1/4,1/2), (1/3,1) circles
% Corners are (5/14,.1237)=(.3571,.1237) ,(7/26,.0666)=(.2692,.0666)
\setlength{\unitlength}{10in} % Note: use one more significant digit
\centerline{ % Center this figure
\begin{picture}(.2,.44)(.2,-.02)
\linethickness{.8pt}
\put(.2,0){\line(1,0){.2}} % x-axis
\put(.375,0){\arc(-.0179,.1237){49.576}} % Three sides of domain
\put(.2667,0){\arc(.0666,0){86.992}}
\put(.6667,0){\arc(-.3096,.1237){21.779}}
\linethickness{.4pt}
\put(.3,0){\line(0,1){.4}} % Vertical geodesic
\put(-.68,0){\arc(1,0){23.578}} % Approximating geodesics
\put(1.28,0){\arc(-.9165,.4){23.578}}
\put(.3,.3162){\vector(0,-1){0}} % Arrowheads on geodesics
\put(.2687,.3162){\vector(1,-3){0}}
\put(.3313,.3162){\vector(-1,-3){0}}
\put(.28,-.02){\clap{$\frac{p_2}{q_2}$}} % Endpoints of geodesics
\put(.3,-.02){\clap{$\theta$}}
\put(.32,-.02){\clap{$\frac{p_1}{q_1}$}}
\put(.3,.41){\vclap{\clap{$\gamma$}}} % Labels on geodesics
\put(.23,.41){\vclap{\clap{$\gamma'_1$}}}
\put(.37,.41){\vclap{\clap{$\gamma'_2$}}}
\put(.33,.1){\rlap{$\sF_j$}} % Label in domain
\put(.3140,.1091){\circle*{.004}} % Points z_{1,j},z_{2,j}
\put(.2836,.0853){\circle*{.004}}
\put(.3340,.1491){\vector(-1,-2){.018}}
\put(.3340,.1491){\rlap{$z_{1,j}$}}
\put(.282,.09){\llap{$z_{2,j}$}}
\end{picture}
}
\caption{Intersecting geodesics: $\gamma$ is trapped between $\gamma'_1$
and $\gamma'_2$ within $\sF_j$.}
\end{figure}
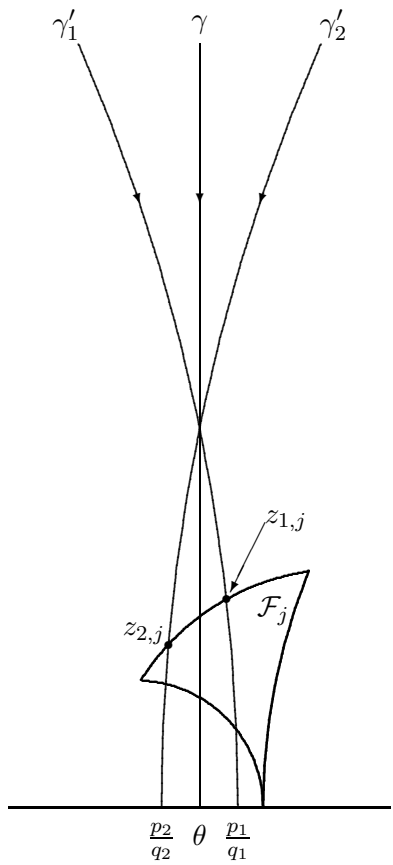
\else %\iffigures
$$
\begin{array}{c}
\hline
~~~~~~~~~~~\mbox{Figure~6.1 about here}~~~~~~~~~ \\
\hline
\end{array}
$$
\fi

There are two main types of forbidden blocks (including initial blocks).
Some blocks are {\em edge-forbidden}; they cannot occur in a cutting
sequence because they correspond to a geodesic hitting two edges which
are part of the same hyperbolic line.  Others are {\em whole-forbidden}
because a geodesic can hit any pair of edges in the sequence, but no
geodesic can hit the whole sequence.

Series~\cite[Theorem 3.1]{Ser86} shows that any minimal edge-forbidden
block for any geodesic flow must correspond to geodesic hitting a
sequence of domains which are all adjacent to the same boundary line,
crossing that line as the first and last step.
This allows us to find the edge-forbidden blocks from Figure~2.3;
without loss of generality, we can let the boundary line be the unit
circle, with the forbidden block starting inside the unit circle, going
outside, and coming back inside.  Since the unit circle has three
segments, there are nine such blocks: $\bar{\hc{J}}\bar{\hc{J}}$,
$\bar{\hc{L}}\bar{\hc{R}}$, $\bar{\hc{R}}\bar{\hc{L}}$,
$\bar{\hc{L}}\bar{\hc{J}}\bar{\hc{L}}\bar{\hc{J}}$,
$\bar{\hc{R}}\bar{\hc{J}}\bar{\hc{R}}\bar{\hc{J}}$,
$\bar{\hc{J}}\bar{\hc{L}}\bar{\hc{J}}\bar{\hc{L}}$,
$\bar{\hc{J}}\bar{\hc{R}}\bar{\hc{J}}\bar{\hc{R}}$,
$\bar{\hc{L}}\bar{\hc{J}}\bar{\hc{L}}\bar{\hc{L}}\bar{\hc{J}}\bar{\hc{L}}$,
$\bar{\hc{R}}\bar{\hc{J}}\bar{\hc{R}}\bar{\hc{R}}\bar{\hc{J}}\bar{\hc{R}}$.
These also lead to excluded initial blocks.  Since the edges
$\bar{\hc{R}}$ and $\bar{\hc{L}}$ meet at infinity, any forbidden block
beginning with either one corresponds to an excluded initial block.
These blocks are exactly those which are forbidden by
Theorem~\ref{CStoOCF}.

We next obtain from Theorem~\ref{1h1mthm} a characterization of the the
whole-forbidden blocks and whole-excluded initial blocks in $\Pi_\sF^0$.
These conditions involve a {\em critical symbol} $1_h$ or $1_m$ which is
mislabeled in a cutting sequence expansion.

\begin{defn}\label{ambigdefn}
A finite sequence $[d_1  \dd  d_n , 1, b_1  \dd  b_m ]$ of positive
integers is {\em ambiguous}
if it does not suffice to determine the MGCF-label on the critical
symbol 1.  More precisely, if
\begin{eqnarray*}
\delta_0  & : = & [0, d_n  \dd  d_1 ], ~ \mbox{and ~$
\delta_1  : =  [0, d_n  \dd  d_1 + 1 ]$},  \\
\beta_0  & : = & [ 1 ,  b_1  \dd b_m ], ~\mbox{and ~ $
\beta_1  : =   [1, b_1  \dd b_m + 1 ]~$ },
\end{eqnarray*}
then $[d_1  \dd  d_n , 1, b_1  \dd b_n ]$ is
{\em ambiguous}
if and only if
\BE
 \interval [ \beta_0  ,  \beta_1 ]  \cap \interval [ \hc{N}(\delta_0)
 ,\hc{N}( \delta_1 )] \neq  \emptyset,
\EE
Here $ \interval [ \delta_0  ,  \delta_1 ]$ denotes the closed interval
determined by $\delta_0$ and
$\delta_1$ with any ordering of the endpoints, i.e.
$\delta_1  <  \delta_0$ may occur.
\end{defn}

Note that the interval $\interval [\hc{N}(\delta_0),\hc{N}( \delta_1 )]$
contains all real numbers whose continued-fraction expansions begin with
the same partial quotients as $\delta_0$.

\begin{defn}
A {\em central sequence} $[d_1,\ldots,d_n,1,b_1,\ldots,b_n]$ is an
ambiguous sequence such that $p/q=[0{\rm\ or\ -1}, d_1 \dd d_n , 1, b_1
\dd b_m ]$ is a
$1_c$ sequence.  We say that $p/q$ is the rational number {\em
associated to} this central sequence.
\end{defn}

We use the following definition to designate initial words of
$\Pi_\sF^0$ using the initial symbol $d_1=\infty$ to indicate that there
is no preceding term in the cutting sequence but that the geodesic comes
from $\infty$.

\begin{defn}
An {\em initial sequence}
$[\infty,d_2,\ldots,d_n,1,b_1,\ldots,b_m]$ denotes an initial word,
which if $-\frac{1}{2}<\theta<0$ has $d_2=1$, $d_3=a_1-1$, and
$d_i=a_{i-2}$ for $i\ge 4$, while if $0<\theta<\frac{1}{2}$,
$d_i=a_{i-1}$ for all $i\ge 2$.
\end{defn}

We define ambiguous initial sequences as in
Definition~\ref{ambigdefn}.  If we set
\BE
\delta_0=\delta_1=[0,d_n,\ldots,d_2]=[0,d_n,\ldots,d_2,\infty],
\EE
then the initial sequence is {\em ambiguous} if
\BE
 \interval [ \beta_0  ,  \beta_1 ]  \cap \interval [ \hc{N}(\delta_0)
 ,\hc{N}( \delta_1 )] \neq  \emptyset,
\EE
where $\interval
[\hc{N}(\delta_0),\hc{N}(\delta_1)]$ is the single point
$\hc{N}(\delta_0)$.

\begin{thm}[Characterization of $\Pi_\sF^0$]\label{blockthm}
The set $\Pi_\sF^0$ of cutting sequences of irrational $\theta$ in
$-\frac{1}{2} < \theta < \frac{1}{2}$ uses the alphabet $\{\bar{\hc{R}},
\bar{\hc{L}}, \bar{\hc{J}} \}$.  It consists of all sequences that
factorize in segments $\hc{B}_0 \hc{B}_1 \hc{B}_2 \ldots$ as in
Theorem~\ref{CStoOCF}, with the additional property that no sequence of
consecutive segments corresponds to a Minkowski geodesic continued
fraction expansion
\[
[d_1 , \ldots , d_n , 1_\ast , b_1 , \ldots , b_m ]
\]
with $d_1=\infty$ allowed, with associated
$\delta_0,\delta_1,\beta_0,\beta_1$, such that:

(i) Both of $[d_2 , \ldots , d_n , 1, b_1 , \ldots , b_m ]$ and
$[d_1 , d_2 , \ldots , d_n , 1, b_1 , \ldots , b_{m-1}]$ are ambiguous.

(ii) Either
\BE
1_\ast = 1_h {\rm\ and\ } \interval [\beta_0 , \beta_1 ] <
\interval [ \hc{N}(\delta_0) , \hc{N} (\delta_1) ],
\EE

or
\BE
1_\ast = 1_m {\rm\ and\ } \interval [\beta_0 , \beta_1 ] >
\interval [ \hc{N}(\delta_0) , \hc{N} (\delta_1) ] ~.
\EE
\end{thm}

{\bf Proof.}
Immediate from Theorem \ref{1h1mthm} using Theorems \ref{MGCFtoOCF}
and~\ref{CStoOCF}.  $\qed$

Theorem~\ref{blockthm} characterizes a large set of forbidden blocks
which are given by conditions (i) and (ii) when $d_1$ is finite.  It
also gives an additional set of excluded initial blocks when
$d_1=\infty$.  These are sufficient to determine $\Pi_\sF^0$ because
they identify each 1 as $1_h$ or $1_m$ as appropriate.  However,
Theorem~\ref{blockthm} does not give a complete set of minimal forbidden
blocks, because there are extra forbidden blocks in which several
symbols $d_i=1$ and $b_i=1$ are replaced by $1_h$ and $1_m$, which if
left as indeterminate $1$'s would not be forbidden.

The complete set of minimal forbidden blocks of $\Pi_\sF^0$ seems harder
to characterize.  Consider a given block of symbols
$\bar{\hc{L}},\bar{\hc{R}},\bar{\hc{J}}$; this can be parsed by
Theorem~\ref{CStoOCF} into a block of complete segments
$[a_1,\ldots,a_m]$ in which each symbol $a_i$ is either $1_h$, $1_m$, or
an integer at least 2, together with some conditions on the adjacent
incomplete segments; for example, an incomplete segment $\bar{\hc{L}}^4$ can be
part of a complete segment encoding either $a_{m+1}\ge 4$, or
$a_{m+1}=3$, $a_{m+2}=1_m$.  Then $[a_1,\ldots,a_m]$ is a forbidden
block if and only if a certain finite set of linear fractional
conditions on two real numbers $\alpha$, $\beta$ of the form
\BE
\alpha<\frac{a\beta+b}{c\beta+d}
\label{consistcond}
\EE
with $ad-bc=\pm 1$ are inconsistent.  There is one such inequality
for each symbol $1_h$ or $1_m$ in the block, which encodes the condition
that $[\alpha,a_1,\ldots,a_m,\beta]$ produces the correct symbol $1_h$
or $1_m$.  There may also be one or two inequalities on $\alpha$ alone
(or on $\beta$ alone) if there is an incomplete segment at that end, and
one more condition~\numb{consistcond} if the possible encoding of that
segment includes another symbol $1_h$ or $1_m$.  The total number of
inequalities is linear in $m$.

Tables 6.1, 6.2, and 6.3 below list some central sequences, ambiguous
sequences and forbidden blocks, respectively; these were obtained by
applying the transformation $\hc{N}$ to simple $\delta_i$.  Table 6.3
illustrates the computation of forbidden sequences from a central
sequence; any change to the terms in a central sequence which increases
either $\alpha_n$ (and thus decreases $\hc{N}(\alpha_n)$) or $\beta_n$
forces the 1 to be $1_h$, and conversely for $1_m$.

\iftables
% \dblfin{3j+3}{1_c,1_h,j,1_h,2}$\ (j<4)$\\ % Too long for table
% \dblfin{3j+3}{1_c,1_h,j,1_m,2}$\ (j\ge 4)$\\
\begin{table}
% \begin{tabular}{|r@{}l|r@{}l|r@{}l|r@{}l|} % In full-width version
\begin{tabular}{|r@{}l|r@{}l|} % In narrow version
\hline
\shortrow{2}{1,1}{4}{3,1}
\shortrow{3}{1,2}{2,2}{2,1,1}
\shortrow{4}{1,3}{2}{1,1}
\shortrow{3j+1}{1,3j}{1,j}{1,j-1,1}
\shortrow{3j+2}{1,3j+1}{1,j,3}{1,j,2,1}
\shortrow{2,2}{1,1,2}{3}{2,1}
\shortrow{2,3}{1,1,3}{2,5}{2,4,1}
\shortrow{3,2}{1,2,2}{3,4}{3,3,1}
\shortrow{4,3}{1,3,3}{2,3}{2,2,1}
\shortrow{5,2}{1,4,2}{3,2}{3,1,1}
\shortrow{j,1}{1,j-1,1}{3j+1}{3j,1}
\shortrow{3,j,1}{1,2,j,1}{3j+2}{3j+1,1}
% \hline Omit in narrow version
\end{tabular}
\caption{Some central sequences, including all with at most five terms.}
\end{table}

\begin{table}
\begin{tabular}{|r@{}l|r@{}l|}
\hline
\symrow j 1 {\un} {3j+1}{{\rm any}}
\symrow j 1 {\un} {3j+2}{{\rm any}}
\symrow j 1 {\un} {3j+3}{{\rm any}}
\dblunf{1,2}{2,1} & &\\
\symrow 2 2 {\un} 3 {\ge 4}
\symrow 3 2 {\un} 3 2
\symrow 3 2 {\un} 3 3
\symrow 4 2 {\un} 3 2
\dblunf{\ge 5,2}{3,1}&\dblunf{1,3}{2,\ge 5} \\
\hline
\end{tabular}
\caption{Non-central ambiguous sequences which go two terms forward and
two terms back from the underlined 1.}
\end{table}

\begin{table}
\begin{tabular}{|r@{\,}c@{\,}l|}
\hline
$[\ge 3,$&$1_h,$&$ \ge 3]$\\
$[\ge 3,$&$1_h,$&$2,1]$\\
\hline
$[1,$&$1_m,$&$ 1]$\\
$[\ge 1,2,$&$1_m,$&$ 1]$\\
$[1,$&$1_m,$&$2, \ge 2]$\\
$[\ge 1,2,$&$1_m,$&$2,\ge 2]$\\
\hline
\end{tabular}
\caption{Forbidden sequences obtained from the central sequence
$[1,2,1_c,2,2]$; the reverses of these sequences are also forbidden.}
\end{table}
\else %\iftables
$$
\begin{array}{c}
\hline
~~~~~~~~~~~\mbox{Tables 6.1, 6.2, and 6.3 about here}~~~~~~~~~ \\
\hline
\end{array}
$$
\fi

Theorem~\ref{blockthm} implies that $\Pi_\sF^0$ is complicated in the
sense that its set of minimal forbidden blocks is very large.

\begin{thm}\label{expthm}
The number $n(k)$ of minimal forbidden blocks of $\Pi_\sF^0$ of
length at most $k$ grows exponentially in $k$; in fact
\BE
\liminf_{k\to\infty}n(k)^{1/k}\ge 2^{1/12}.
\label{explimit}
\EE
\end{thm}

{\bf Proof.} We will show that each central sequence associated to a
rational $p/q\ne 1/2$ yields two minimal forbidden blocks.  The
forbidden blocks are produced by adding one symbol to each end of the
central sequence, and by replacing the central $1_c$ with a three-symbol
word.
All these forbidden blocks are distinct.  Assuming these facts are proved,
consider the
central sequences
\BE
[d_1,d_2,\ldots,d_n,1_c,b_1,\ldots,b_m].
\label{expcenter}
\EE
in which each $d_i=1$ or 2, and in which $[b_1,\ldots,b_m]$ is
determined from $[d_1,\ldots,d_n]$.  The number of symbols in the cutting
sequence encoding of $[d_1,\ldots,d_n]$ is at most $3n$, and by
Corollary~\ref{autfinds1c}, the number of symbols in the cutting
sequence encoding of $[b_1,\ldots,b_m]$ is at most
$9n$.  In obtaining forbidden blocks, the $1_c$ term is encoded by three
symbols, and one symbol is added at each end, hence all resulting
forbidden blocks contain at
most $12n+5$ symbols.  We conclude that there are at least
$2^{n+1}$ minimal forbidden blocks of length at most $12n+5$, which
proves~\numb{explimit}.

We now construct the minimal forbidden blocks.  We use the result of
Appendix~A, which shows that any vertical geodesic for $\theta=p/q$ with
$-\frac{1}{2}<\theta<\frac{1}{2}$ hits at most one corner of a fundamental
domain.  This
fact implies that each of the central sequences~\numb{expcenter}
contains only one symbol $1_c$, and thus each $p/q$ comes from at most
one central sequence.
  Associated to the central sequence is a word
$\hc{W}_1\hc{S}\hc{W}_2$ in which $\hc{S}=\bar{\hc{C}}_1$ or
$\bar{\hc{C}}_2$ is a corner symbol.  For the corner symbol, there are
two choices of three symbols in
$\{\bar{\hc{R}},\bar{\hc{L}},\bar{\hc{J}}\}$, which replace the $1_c$ by
$1_h$ or $1_m$; for $\bar{\hc{C}}_1$, the choices are
$\bar{\hc{J}}\bar{\hc{R}}\bar{\hc{J}}$ or
$\bar{\hc{L}}\bar{\hc{J}}\bar{\hc{L}}$.  We consider the eight words
obtained from $\hc{W}_1\hc{S}\hc{W}_2$ by replacing $\hc{S}$ by either
choice of a three-symbol block, and by adding a prefix symbol and a
suffix symbol, each of which may be either $\bar{\hc{L}}$ or
$\bar{\hc{R}}$.  The continued fraction for $p/q$ can be recovered from
any one of these eight blocks, and since it has only a single $1_c$, the
critical $1_c$ in the central sequence is also uniquely determined.

We claim that six of these eight blocks are admissible blocks for
$\Pi_\sF^0$ and the other two are forbidden blocks.  The vertical
geodesic $\gamma=\langle\infty,p/q\rangle$ corresponding to
\[
\frac{p}{q}=[0{\rm\ or\ }-1,d_1,\ldots,d_m,1_c,b_1,\ldots,b_m]
\]
can be approximated by geodesics in the symbol topology in six different
ways.  Two of these consist of approximating geodesics which do not
cross $\gamma$ at all but approach it from the left and right,
respectively.  The other four consist of approximating geodesics which
cross $\gamma$, either crossing above or below the corner that $\gamma$
hits at the $1_c$, and initially approaching $\gamma$ either from the
left or from the right.  For sufficiently good approximations, these
produce the admissible blocks.  (Here we again use the fact that the geodesic
$\gamma$ hits exactly one corner, which implies that all geodesics
sufficiently close to $\gamma$ have the same sequences $\hc{W}_1$ and
$\hc{W}_2$.)
 The other two blocks are forbidden by condition (ii) of
Theorem~\ref{blockthm}.  They encode cutting sequences for a geodesic
that would have to approach $\gamma$ from the left, pass to the right of
the corner, and then return to the left of $\gamma$, or vice versa; such
a geodesic would have to cross $\gamma$ twice, which is impossible.

Let $\bar{\hc{C}}_a$ and $\bar{\hc{C}}_b$ denote the two possible three-symbol
encodings of $\hc{S}$.  Of the four blocks
$$
\bar{\hc{R}}\hc{W}_1\bar{\hc{C}}_{a}\hc{W}_2\bar{\hc{R}}~,~~
\bar{\hc{R}}\hc{W}_1\bar{\hc{C}}_{a}\hc{W}_2\bar{\hc{L}}~,~~
\bar{\hc{L}}\hc{W}_1\bar{\hc{C}}_{a}\hc{W}_2\bar{\hc{R}}
~~\mbox{and}~~
\bar{\hc{L}}\hc{W}_1\bar{\hc{C}}_{a}\hc{W}_2\bar{\hc{L}}~,
$$
exactly three are admissible and one is forbidden.  Every sub-block of
the forbidden block appears in one of the three admissible blocks, hence
the forbidden block is minimal.  The same argument applies to the other
four blocks containing $\bar{\hc{C}}_b$, and produces a second minimal
forbidden block.  $\qed$

\begin{thm}\label{noautthm}
There does not exist a finite automaton which, when given the additive
continued fraction expansion of an irrational number $\theta$ with $0 <
\theta < \frac{1}{2}$ as its input sequence, computes the Minkowski
geodesic continued fraction expansion of $\theta$, or, equivalently, the
cutting sequence expansion of the geodesic $\langle \In , \theta
\rangle$.
\end{thm}

{\bf Proof.}
Any finite automaton that would compute the Minkowski geodesic expansion
of $\theta$ must
output the $n$-th term of this expansion after seeing at most a bounded
number of symbols following the $n$-th symbol of the ACF expansion of $\theta$.
However the forbidden block criteria of Theorem~\ref{blockthm} show that
it is sometimes
necessary to see an arbitrarily large string of symbols after the
$n$-th symbol, to decide if $1_h$ or $1_m$ should be used.
For example, for each $j\ge 1$, the sequence
\BE
[2^{4j+2},1_c,3,(8,4)^j] \label{badfin}
\EE
is a central sequence.  Adding some later terms decreases
$\beta_{4j+2}$, and thus changes the $1_c$ to $1_m$,  while decreasing
the final 4 in $(8,4)^j$  to a 3 and then adding some further terms
increases
$\beta_{4j+2}$ and thus changes the $1_c$ to $1_h$.  A finite automaton
cannot look ahead through the $14j+6$ terms which are necessary, since
$j$ can be any integer.  In
fact, even a pushdown automaton (with one stack) cannot correctly
compute all such cutting sequences, because it must look $14j+6$ steps
ahead, but also
$12j+6$ steps back to distinguish $[a_0;2^{4j},1_m,3,(8,4)^j,8,3,\ldots]$
from $[a_0;2^{4j+2},1_h,3,(8,4)^j,8,3,\ldots]$.
$\qed$

\section{Two-sided Cutting Sequences: Structure of $\Sigma_\sF$}
\hsp
We now use the information on vertical cutting sequences $\Pi_\sF^0$
to characterize the two-sided cutting sequences $\Sigma_\sF$.

\begin{thm}\label{twosidethm}
The cutting sequence shift $\Sigma_\sF$ for the fundamental domain
$\sF$ of $PSL (2, \ZZ )$ is the closed subshift whose forbidden blocks
coincide with the set of forbidden blocks of $\Pi_\sF^0$.
\end{thm}

{\bf Proof.} This result follows from the fact that the set of images of
the set of irrational vertical geodesics under $PSL(2,\ZZ)$ is dense in
the space of all geodesics of $\hpl /PSL (2, \ZZ )$.

Let $\sS_\Pi^0$ and $\sS_\Sigma^0$ and $\sS_\Sigma$ denote the complete
sets of forbidden blocks of $\Pi_\sF^0$, $\Sigma_\sF^0$, and
$\Sigma_\sF$, respectively.  Since $\Sigma_\sF$ is the closure of
$\Sigma_\sF^0$, we have $\sS_\Sigma^0=\sS_\Sigma$.

We first show that
\BE
\sS_\Sigma^0 \subseteq \sS_\Pi^0.
\EE
For this it suffices to show that every admissible word $\hc{W}$ in
$\Pi_\sF^0$ is an admissible word in $\Sigma_\sF^0$.  Suppose that the
word $\hc{W}$ appears in the cutting sequence of the irrational vertical
geodesic $\langle\infty,\theta\rangle=\{\theta+it:t>0\}$.  This geodesic
is a limit of geodesics $\langle\phi_i,\theta\rangle$ where
$\phi_i\to\infty$ through a sequence of values such that
$\langle\phi_i,\theta\rangle$ hits no corner of an image of $\sF$.  The
word $\hc{W}=\hc{S}_1\ldots\hc{S}_r$ corresponds to a specific set of
edges of translated fundamental domains $\{g_j\sF:0\le j\le r\}$ with
$g_j\in SL(2,\ZZ)$.  For all sufficiently large $\phi_i$, the geodesic
$\langle\phi_i,\theta\rangle$ passes through the same sequence of
fundamental domains $\{g_j\sF:0\le j\le r\}$, hitting the same sequence
of edges in the same order.  Thus $\hc{W}$ occurs in the two-sided
cutting sequence of $\langle\phi_i,\theta\rangle$, so it is an
admissible word of $\Sigma_\sF^0$.

The reverse inclusion
\BE
\sS_\Pi^0 \subseteq \sS_\Sigma^0. \label{revincl}
\EE
is proved similarly.  Let $\hc{W}$ be an admissible word in some general
position geodesic $\gamma=\langle\phi,\theta\rangle$.  There is a
sequence of translated fundamental domains $\{\sF_j:0\le j\le r\}$
which $\phi$
passes through that corresponds to $\hc{W}$.  Since $\gamma$ hits no
corners, we can choose a rational number $p/q$ sufficiently close to
$\phi$ such that the geodesic $\gamma'=\langle p/q,\theta\rangle$ passes
through the same set of fundamental domains $\{\sF_j:0\le j\le r\}$,
hitting the same sequence of edges in the same order, hence its cutting
sequence (which is only one-sided infinite) contains the word $\hc{W}$.
Now apply to $\gamma'$ the transformation $\hc{M}$ in $PSL(2,\ZZ)$ which
takes $p/q$ to $\infty$, and $\theta$ to some $\theta'$ in
$(-\frac{1}{2},\frac{1}{2})$.  For this choice, the cutting sequence of
$\hc{M}(\gamma')$
is in $\Pi_\sF^0$, and~\numb{revincl} follows.  $\qed$

%Generalized continued fraction removed here

\begin{thm}\label{encodingthm}
Every element in $\Sigma_\sF$ is a cutting sequence for a unique
oriented geodesic on $\hpl /PSL (2, \ZZ )$ which hits $\sF$, except for
the two sequences $\bar{\hc{R}}^\infty$ or $\bar{\hc{L}}^\infty$.  Every
oriented geodesic $\gamma$ on $\hpl /PSL(2, \ZZ)$ has at least one and at
most finitely many shift-equivalence classes of cutting sequences in
$\Sigma_\sF$.
If $\gamma$ is not periodic then it has at most eight shift-equivalence
classes of cutting sequences in $\Sigma_\sF$.
\end{thm}

To establish this result, we first prove a preliminary lemma.

\begin{lemma}
Let $\gamma_j=\langle\theta'_j,\theta_j\rangle$ for $j=1,2,\ldots$ be a
sequence of general position geodesics that intersect $\sF$ which have
cutting sequences $\{\hc{S}_i^{(j)}:i\in\ZZ\}$ such that the symbol
$\hc{S}_0^{(j)}$ corresponds to the geodesic $\gamma_j$ entering the
fundamental domain $\sF$.  If the cutting sequences $\{\hc{S}_i^{(j)}\}$
converge in the symbol topology as $j\to\infty$ to a limit sequence
$\{\hc{S}_i\}$ then the endpoints $\theta'_j$ and $\theta_j$
converge to unequal limiting values
\BE
\theta'=\lim_{j\to\infty}\theta'_j{\rm\ and\ }
\theta=\lim_{j\to\infty}\theta_j.
\label{geodlimits}
\EE
in $\RR\cup\{-\infty,\infty\}$.
The geodesics $\gamma_j$ converge to a limiting geodesic
$\gamma=\langle\theta',\theta\rangle$ if at least one of $\theta,\theta'$ is
finite.  The exceptional cases $\langle-\infty,\infty\rangle$ and
$\langle\infty,-\infty\rangle$ correspond to the limiting symbol
sequences $\bar{\hc{R}}^\infty$ and $\bar{\hc{L}}^\infty$, respectively.
\end{lemma}

{\bf Proof.}  Let $\hc{S}_+=\{\hc{S}_i:i>0\}$ and
$\hc{S}_-=\{\hc{S}_i:i\le 0\}$ denote the
positive and negative part of the limit sequence, respectively.  If
$\hc{S}_+$ is $\bar{\hc{R}}^\infty$, then since any geodesic
$\gamma_j$ which enters $\sF$ after $\hc{S}_0$ and has
$\hc{S}_i^{(j)}=\bar{\hc{R}}$ for $1\le i\le n$ must have $\theta_j\ge
n+1/2$, we conclude that $\lim_{j\to\infty}\theta_j=\infty$ in this
case.  Similarly, if $\hc{S}_+$ is $\bar{\hc{L}}^\infty$,
then $\lim_{j\to\infty}\theta_j=-\infty$; if $\hc{S}_-$ is
$\bar{\hc{R}}^\infty$, then $\lim_{j\to\infty}\theta'_j=-\infty$ if
$\hc{S}_-$ is $\bar{\hc{L}}^\infty$, then
$\lim_{j\to\infty}\theta'_j=-\infty$.  In particular, this shows that
the two-sided limit sequence $\bar{\hc{R}}^\infty$ corresponds to
$\langle-\infty,\infty\rangle$, and $\bar{\hc{L}}^\infty$ corresponds to
$\langle\infty,-\infty\rangle$.

We next observe that since $\bar{\hc{R}}\bar{\hc{L}}$ and
$\bar{\hc{L}}\bar{\hc{R}}$ are forbidden blocks, if $\hc{S}_+$ is not
$\bar{\hc{R}}^\infty$ or $\bar{\hc{L}}^\infty$ it must contain a symbol
$\bar{\hc{J}}$.  Let the initial segment of $\hc{S}_+$ up to
the first such symbol be $\bar{\hc{R}}^n\bar{\hc{J}}$ (resp.
$\bar{\hc{L}}^n\bar{\hc{J}}$).  Any geodesic $\gamma_j$ which matches
these symbols necessarily has $\theta_j\le n+1$ since it crosses the
semicircle with endpoints $n-1$ and $n+1$ (resp. $\theta_j\ge
-n-1$) and is positively oriented (resp. negatively oriented), hence
$\theta_j\ge-1/2$ (resp. $\theta_j\le 1/2$).  In either case, if a
limit~\numb{geodlimits} exists it must be finite.  Similar reasoning
applies to $\hc{S}_i$ to show that
$\lim_{j\to\infty}\theta'_j$ is infinite if and only if
$\hc{S}_-$ does not contain a symbol $\bar{\hc{J}}$.

We now suppose that $\hc{S}_+$ contains a symbol $\bar{\hc{J}}$, and
claim that $\lim_{j\to\infty}\theta_j$ exists and is finite.  The
argument above shows that $\{\theta_j\}$ is bounded, so to prove this
claim it suffices to show that if $\{\theta_j\}$ is not a Cauchy
sequence then the one-sided sequences $\{\hc{S}_i^{(j)}:i> 0\}$ cannot
converge in the symbol topology.  If it is not a Cauchy sequence, there
is some $\epsilon>0$ such that for any $N$, we have $j,k\ge N$ with
$|\theta_j-\theta_k|>\epsilon$.  The general-position geodesics
$\gamma_j$ and $\gamma_k$ enter $\sF$ and thus have radius at least
$\sqrt{3}/2$.  Consequently, any point $z=x+iy$ on the geodesic
$\gamma_j$ with $|x-\theta_j|<1/4$ and $0<y<\epsilon'$ for
$\epsilon'<\epsilon^2/36$ will actually have $|x-\theta_j|<\epsilon/3$,
and thus $|x-\theta_k|>2\epsilon/3$; similarly, points on $\gamma_k$
with $0<y<\epsilon'$ and $|x-\theta_k|<1/4$ must have
$|x-\theta_k|>2\epsilon/3$.  This implies an upper bound $N_\epsilon$ on
the number of fundamental domains which $\gamma_j$ and $\gamma_k$ both
hit somewhere inside the box $|Re(z)-(\theta_j+\theta_k)/2|\le
2\epsilon$, $Im(z)\le \epsilon$.  We may assume $\epsilon<1/2$, hence
any such fundamental domain has a cusp at a finite rational point $p/q$.
The rosette $R(p/q)$ of all translated fundamental domains that touch
$p/q$ has Euclidean diameter less than $4/q^2$.  Thus any translated
fundamental domain $\sF'$ which is hit by both $\gamma_j$ and $\gamma_k$
inside the box must have two points whose real parts differ by at least
$\epsilon/3$, hence $q<\sqrt{12/\epsilon}$.  In addition, $p/q$ must lie
within Euclidean distance $4/q^2$ of both $\theta_j$ and $\theta_k$,
hence there are only finitely many such rational $p/q$.  Finally, there
is a finite bound depending only on $q$ and $\epsilon$ on the number of
fundamental domains in each rosette $R(p/q)$ which have Euclidean
diameter exceeding $\epsilon/3$.  Together this yields a finite upper
bound $N_\epsilon$ on the number of common fundamental domains $\sF'$
between $\gamma_j$ and $\gamma_k$ in the box.  Outside the box, after
exiting $\sF$ the geodesics $\gamma_j$ and $\gamma_k$ traverse the
strip $\epsilon'\le Im(z)\le 1$, and thus can hit at most
$N'_\epsilon$ translated fundamental domains $\sF'$ in that strip,
because there are only finitely many geodesics with Euclidean radius at
least $\epsilon'$ which intersect the region
$-\frac{1}{2}<Re(z)<\frac{1}{2}$.  We
conclude that the positive cutting sequences $\gamma_j$ and $\gamma_k$
cannot agree on all of their first $N_\epsilon+N'_\epsilon+n+1$ symbols.
Thus the sequences $\hc{S}_+^{(j)}=\{\hc{S}_i^{(j)}:j>0\}$ cannot
converge in the symbol topology, which proves the claim.

A similar argument applies to the limit sequence $\hc{S}_-$
if it contains the symbol $\bar{\hc{J}}$, to prove that
$\lim_{j\to\infty}\theta'_j=\theta'$ exists and is finite.

The geodesic $\gamma=\langle\theta',\theta\rangle$ is encoded as a limit
of the sequence of geodesics $\gamma_i$.  If all of the $\gamma_i$ for
$i\ge N$ intersect a particular domain $\sF'$, then since the endpoints
$\theta'_i$ and $\theta_i$ converge to $\theta'$ and $\theta$, and $\sF$
is closed, $\gamma$ also intersects $\sF$, at least on the boundary.  If
no $\gamma_i$ for $i\ge N$ intersect a domain $\sF'$, then $\gamma$ does
not intersect the interior of $\sF'$.  Thus, if we consider $\gamma$ to
hit a domain if either it passes through the interior, or it intersects
the boundary and all $\gamma_i$ for sufficiently large $i$ intersect the
interior, we can define a cutting sequence as the set of edges which
separate these domains, and this sequence is the limit in the symbol
topology of the cutting sequences of $\gamma_i$.  $\qed$

{\bf Proof of Theorem~\ref{encodingthm}}.  Since $\Sigma_\sF$ is a
compact set, every oriented geodesic
$\gamma=\langle\theta',\theta\rangle$ which hits the interior of $\sF$
has at least one cutting sequence in $C(\gamma)$ by taking the limit
point of cutting sequences from a family
$\gamma_j=\langle\theta'_j,\theta_j\rangle$ with $\theta'_j\to\theta'$
and $\theta_j\to\theta$ as $j\to\infty$.

It remains to show that each oriented geodesic in $\hpl$ which hits $\sF$
has only a finite number of shift-equivalence classes of cutting
sequences in $[C(\gamma)]$.  The simplest cases are general position
geodesics, for which $C(\gamma)$ is a single cutting sequence, and
$[C(\gamma)]$ consists of a single shift-equivalence class.

We determine how many shift-equivalence classes of cutting sequences are
possible for a
limiting geodesic.  First, consider the limiting geodesics with exactly
one rational endpoint (considering $\infty$ as rational).  These
correspond to vertical geodesics under the $PSL(2,\ZZ)$ action, so by
Lemma~\ref{vertcorner}, such geodesics cannot hit any internal corners,
and thus have two possible encodings; either can be obtained by using
approximating geodesics which approach the rational endpoint from either
side.  Thus $[C(\gamma)]$ contains two shift-equivalence classes.

Next, consider those geodesics which have two rational endpoints.
Without loss of generality, we can move one endpoint to $\infty$ and so
have a vertical geodesic $\gamma=\langle\infty,p/q\rangle$.  Such a
geodesic hits only finitely many corners of fundamental domains.  If it
hits $n$ corners at finite values of $t$, then $C(\gamma)$ contains
exactly $2n+4$ cutting sequences.  Two of them come from geodesics
approximating $\gamma$ from either side without crossing it, and the
other $2n+2$ result from approximating geodesics which cross $\gamma$
between the $k$th and $(k+1)$st corners with $0\le k\le n$, either from
left to right or from right to left.  (Here the 0th corner is
$\infty$, and the $(n+1)$st corner is $p/q$.)  In Appendix~A, we prove
that $n\le 1$ for all rational $p/q$, except those $p/q\equiv
\frac{1}{2}\pmod{1}$,
which have $n=2$.  Thus $[C(\gamma)]$ contains at most eight
shift-equivalence classes.

It remains to bound the number of shift-equivalence classes of cutting
sequences for
geodesics which have two irrational endpoints.  If such a geodesic hits
no corners, then it has only one cutting sequence in $C(\gamma)$.  If it
hits exactly one corner, then it has exactly two cutting sequences,
which are obtained using geodesics which approach it while staying on
opposite sides of the corner.  The difficult case occurs with geodesics
that hit at least two corners.  We show that all such geodesics hit
infinitely many corners and are periodic.  In this case, $C(\gamma)$
will be an infinite set.  To show periodicity, we observe that a corner
in the upper half-plane is the intersection of two circles with centers
at rational points on the $x$-axis and rational radii, so its
$x$-coordinate is rational and its $y$-coordinate is the square root of
a rational number.  Thus, if the circle with radius $r$ and center
$(x_0,0)$ passes through two such points $(x_1,y_1)$ and $(x_2,y_2)$, then it
satisfies the equations
\begin{eqnarray*}
(x_1-x_0)^2+y_1^2 &=& r^2, \\
(x_2-x_0)^2+y_2^2 &=& r^2.
\end{eqnarray*}
Since $y_1^2$ and $y_2^2$ are rational, so is $r^2$; also, equating the
left sides gives a linear equation for $x_0$ with rational coefficients.
Thus the circle intersects the $x$-axis in two algebraically conjugate
real quadratic surds.  Pell's equation allows us to find
$M\in SL(2,\ZZ)$ which preserves one endpoint (and
thus its conjugate). Applying this transformation to our
geodesic re-scales it while preserving the orientation; hence the
geodesic necessarily is a periodic geodesic on $\hpl/PSL(2,\ZZ)$.  The
set $C(\gamma)$ is then infinite because we can choose a set of
approximating geodesics which cross between any pair of consecutive
corners of $\gamma$, and the resulting limit cutting sequences are all
distinct.

We complete this case by showing that if a periodic geodesic $\gamma$ on
$\hpl/PSL(2,\ZZ)$ hits exactly $n$ corners of translates of fundamental
domains in its period, then there are exactly $2n+2$
shift-equivalence classes in $[C(\gamma)]$.  Label the corners which are
hit in one period $c_1,\ldots,c_n$.  The approximating geodesics
$\gamma_i$ can only be close in the symbol topology if they all cross
$\gamma$ between the same pair of corners, or if none cross $\gamma$ at
all.  If they cross, then we can use the periodicity of $\gamma$'s
cutting sequence to shift the crossing point between $c_1$ and
$c_{n+1}$.  With $n$ possible crossing regions, and crossing possible
either from inside to outside or vice versa, there are $2n$ encodings;
we get two more from approximating geodesics which are completely inside
or outside $\gamma$, for a total of exactly $2n+2$ shift-equivalence classes.
$\qed$

Theorem~\ref{encodingthm} gives no uniform upper bound on the number of
shift-equivalence classes of
cutting sequences that correspond to periodic geodesics. We formulate
the question whether a uniform upper  bound exists as open problem (3) in
the concluding section.
There are periodic geodesics which
correspond to 10 cutting sequence shift-equivalence classes because they
hit four corners in one
period, such as $\langle-\sqrt{13},\sqrt{13}\rangle$ and
$\langle-\sqrt{133},\sqrt{133}\rangle$.  
Could this be the
maximal number that occurs?  

Theorem~\ref{twosidethm} implies that the shift~$\Sigma_\sF$ is a
relatively complicated set.  Indeed Theorem~\ref{expthm} showed that the
set of minimal forbidden blocks is very large.
We now show that $\Sigma_\sF$ is not a sofic shift.  A {\em sofic shift}
is any shift that is a factor of a shift of finite type (see Marcus and
Lind~\cite{Lin95}.)  Alternatively, it is the set of possible
bi-infinite walks on an edge-labeled finite graph.

\begin{thm}\label{notsofic}
The cutting sequence shift $\Sigma_\sF$ for the fundamental domain $\sF$
of $PSL (2, \ZZ )$ is not a sofic shift.
\end{thm}

{\bf Proof.}
If $\Sigma$ is a shift, the {\em follower set}
$F_\Sigma(\hc{W})$ of a one-sided word
$\hc{W}=(\ldots,\hc{S}_{-2},\hc{S}_{-1})$ is the set of all one-sided
words $\hc{W}^+=(\hc{S}_0,\hc{S}_1,\ldots)$ such that $\hc{W}\hc{W}^+$
is an element of $\Sigma$.  Sofic shifts are characterized by the
property that the totality of different follower sets
$\{F_\Sigma(\hc{W}):{\rm all}\ \hc{W}\}$ is finite; see Marcus and
Lind~\cite[Theorem 3.2.10]{Lin95}.

We show that $\Sigma_\sF$ has infinitely many follower sets.  The
sequences $[3,2^{4j+2},1_c,3,(8,4)^j,10]$ and
$[4,2^{4j+2},1_c,(8,4)^j,13]$ are central.  Thus, in
$[\ldots,3,2^{4j+2},1_*,3,(8,4)^k,13,\ldots]$, the $1_*$ is $1_h$ if
$j\ge k$ but $1_m$ if $j<k$, regardless of the symbols on either side.
Thus the follower sets of one-sided words which end $[3,2^{4j+2}]$ are
different for all $j$.  $\qed$

\section{The Shift $\Sigma_\sF$ Determines $\sF$ up to Isometry}
\hsp
Our object is to prove the following result.

\begin{thm}\label{uniquepoly}
Let $\Gamma$ be a finitely generated discrete subgroup of $PSL (2, \RR )$
that acts properly discontinuously on $\hpl$, and which has a polygonal
fundamental
domain $\sP$ which is hyperbolically convex.
Suppose that the cutting sequence shift $\Sigma_\sP$ is
isomorphic to $\Sigma_{\sF, PSL(2, \ZZ )}$ by a
permutation of symbol alphabets.
Then there is an element $g \in PSL (2, \RR )$ such that
$\sF = g \sP$ and
$$
PSL (2, \RR ) = g \Gamma g^{-1} ~.
$$
\end{thm}

{\bf Proof.}
Since the cutting sequence contains only three symbols, the polygon must
be a triangle.

Since a geodesic can hit the $\bar{\hc{L}}$ edge an unlimited number of
consecutive times, the $\bar{\hc{L}}$ and $\bar{\hc{L}}^{-1}$ edges of
the triangle must intersect at an angle of zero.  The
$\bar{\hc{L}}^{-1}$ edge cannot be the same as the $\bar{\hc{L}}$ edge
because a geodesic can hit it twice consecutively, nor can it be the
$\bar{\hc{J}}$ edge because no geodesic can hit that edge twice
consecutively; thus it must be the $\bar{\hc{R}}$ edge.  Since a
geodesic cannot hit the $\bar{\hc{J}}$ edge twice consecutively, the
$\bar{\hc{J}}$ edge must equal the $\bar{\hc{J}}^{-1}$ edge.  Thus the
generator at the $\bar{\hc{J}}$ edge is an involution, and to have
determinant $1$, it must be an inversion.

Let the angle between $\bar{\hc{J}}$ and $\bar{\hc{L}}$ be $\alpha$, and
the angle
between $\bar{\hc{J}}$ and $\bar{\hc{R}}$ be $\beta$.  If
$2\alpha+\beta<\pi$, then
a geodesic approaching the $\hc{JR}$-corner could hit edges
$\bar{\hc{J}},\bar{\hc{L}},\bar{\hc{J}},\bar{\hc{L}}$ in sequence; however,
this cutting
sequence is not possible for our fundamental domain.  Likewise,
$2\beta+\alpha\ge\pi$.  Also, for any hyperbolic triangle, the sum of
the angles is less than $\pi$, so $\alpha+\beta<\pi$.

The conditions for a given polygon to be a fundamental domain are given
in Maskit~\cite[section 2]{Mas71}.  In our case, it is necessary that
either $n(\alpha+\beta)=2\pi$ for some $n$, or $n\alpha=\pi$ and
$m\beta=\pi$ for some $m$ and $n$.  The only cases consistent with the
conditions on $\alpha$ and $\beta$ are $\alpha=\beta=\pi/3$, the desired
fundamental domain; and $\alpha=\pi/3,\beta=\pi/2$, which is the half of
our fundamental domain with $x>0$.  But that domain is only a
fundamental domain for a group including the reflection in the line
$x=0$, which has determinant $-1$.  (That domain also has different
symbolic dynamics; a geodesic can hit $\bar{\hc{L}}$ and $\bar{\hc{R}}$
consecutively.)  $\qed$

\section{Open Problems}
\hsp
(1).
We have shown in one special case that the polygon $\sP$ can be recovered
from the data $\Sigma_\sP$ up to isometry.
The example $\sF$ comes from a Riemann surface of genus 0.
Does the result persist for higher-genus Riemann surfaces?
Can any $\Sigma_\sP$ be explicitly determined for a
Riemann surface of genus at least one?

(2).
Since $\Sigma_\sF$ determines $\sF$ up to isometry, in principle  it
determines $vol ( \hpl / PSL (2, \ZZ ))$.
Can $vol ( \hpl / PSL (2, \ZZ )) = \frac{ \pi}{3}$ be easily computed
directly from $\Sigma_\sF$?

(3).
Is there a universal upper bound on the number of
shift equivalence classes of cutting sequences
corresponding to any periodic geodesic
on  $\hpl / PSL (2, \ZZ)$?  Equivalently,
is there a universal upper bound on the number of times that a periodic
geodesic  can hit a corner of the  fundamental domain $\sF$, 
during a single period?

(4).
The zeta function $\zeta_\Sigma(z)$ of a shift $\Sigma$ is defined by
\[
\zeta_\Sigma(z)=\exp(\sum_{k=1}^\infty N_k \frac{z^k}{k}),
\]
in which $N_k$ counts the number of periodic words in $\Sigma$ of period
$k$.  Is there a simple formula for the (dynamical) zeta function of
$\Sigma_\sF$?  What is the topological entropy of $\Sigma_\sF$?

(5).
Cutting sequence shifts $\Sigma_\sP$ can be constructed in
higher-dimensional cases along the lines considered in \cite{Lag94},
if one restricts to a suitable subclass of geodesics, called ``flat''
geodesics in \cite{Lag94}.
Can any such $\Sigma_\sP$ be determined explicitly?

{\bf Acknowledgment.}
We are indebted to L.~Flatto and M. Sheingorn
for helpful comments and references.
\newpage
\noindent

% Symbolic dynamics definitions removed here

\appendix
\section{A Bound on the Number of Corners on a Vertical Geodesic}

\begin{lemma}
For a rational $\theta$, the vertical geodesic $\gamma=\{\theta+it:
t>0\}$ has at most one value of $t$ such that $\theta+it$ is a corner of
a $PSL(2,\ZZ)$-translate of $\sF$, unless $\theta\equiv\frac{1}{2}
\pmod{1}$, in which case it has exactly two such values, which are
$t=\sqrt{3}/2$ and $\sqrt{3}/6$.
\end{lemma}

{\bf Proof.}  Since the corner
$-1/2+\sqrt{-3}/2$ of $\sF$ is obtained from the corner
$1/2+\sqrt{-3}/2$ by the
transformation $z\to z-1$, every corner is a $PSL(2,\ZZ)$-translate of
$1/2+\sqrt{-3}/2$.  Let the element of $SL(2,\ZZ)$ be
$\twobytwo{a}{b}{c}{d}$, with $ad-bc=1$.  We then have
\begin{eqnarray}
\twobytwo{a}{b}{c}{d}\left(\frac{1}{2}+\frac{\sqrt{-3}}{2}\right)
&=& \frac{a\left(\frac{1}{2}+\frac{\sqrt{-3}}{2}\right)+b}
         {c\left(\frac{1}{2}+\frac{\sqrt{-3}}{2}\right)+d} \nonumber\\
&=& \frac{(4ac+2ad+2bc+4bd)+(2ad-2bc)\sqrt{-3}}
         {4(c^2+cd+d^2)}\nonumber\\
&=& \frac{(2ac+ad+bc+2bd)+\sqrt{-3}}
	 {2(c^2+cd+d^2)}.\label{corner}
\end{eqnarray}

Let $N=2ac+ad+bc+2bd$ and $D=c^2+cd+d^2$, so that the real part
of~\numb{corner} is $N/2D$.  We will show that $N/2D$ is either in
lowest terms or can be reduced to lowest terms by dividing both $N$ and
$D$ by 3.  (The factor 3 can occur; for example, take $a=2$, $b=1$,
$c=1$, $d=1$.)  Since $ad-bc=1$, $ad+bc$ is odd and thus $N$ is not
divisible by 2.  Now substitute $a=(1+bc)/d$ in~\numb{corner}; this
gives
\begin{eqnarray}
\frac{N}{2D}
&=&\frac{2\frac{1+bc}{d}c+2bd+\frac{1+bc}{d}d+bc}{2D} \nonumber\\
&=&\frac{2c+2bc^2+2bd^2+d+2bcd}{2d(c^2+cd+d^2)} \nonumber\\
&=&\frac{b(c^2+cd+d^2)+2c+d}{2d(c^2+cd+d^2)}.\label{expandcorner}
\end{eqnarray}
Since $(c,d)=1$, we have $c^2+cd+d^2\equiv 1\pmod{2}$, and also
$(d,c^2+cd+d^2)=(d,c^2)=1$.  It follows that
\begin{eqnarray*}
&(N,D)=(dN,D)=(bD+2c+d,D)=(2c+D,D)=(2c+d,c^2+cd+d^2)\\
&=(2c+d,c^2+cd+d^2-d(2c+d))=(2c+d,c(c-d)).
\end{eqnarray*}
Now, $(c,2c+d)=(c,d)=1$, and $(2c+d,c-d)=(3c,c-d)\le(3,c-d)(c,c-d)=1$ or
3.  Therefore $(N,D)=1$ or 3.

Suppose the geodesic is $\theta+it$ with $\theta=v/w$.  There are only
two possible values of~\numb{corner} on this geodesic; $N/2D$ can only
equal $v/w$ in lowest terms if $w=2D$ or $w=2D/3$.  The imaginary part
of~\numb{corner} is either $\sqrt{-3}/w$ or $\sqrt{-3}/3w$.  However,
these two points are separated by a hyperbolic distance of $\ln{3}$, the
length of the finite side of the domain $\sF$; therefore, they can both
be $PSL(2,\ZZ)$-translates of corners only if they are connected by a
$PSL(2,\ZZ)$-translate of an edge.  This edge is the geodesic
$\theta+it$ itself, and the only vertical edges of translates $g\sF$
occur when $\theta=(n+1/2)$ for $n\in\ZZ$.

For $\theta=n+1/2$, there are two $SL(2,\ZZ)$-translates of corners.
If we take $a=1$, $b=n$, $c=0$, $d=1$, we get the corner
$(n+1/2)+\sqrt{-3}/2$; if we take $a=n+1$, $b=n$, $c=1$, $d=1$, we get the
corner $(n+1/2)+\sqrt{-3}/6$.  $\qed$

\begin{tabular}{ll}
{\bf email:} & {\tt jcl@research.att.com} \\
 & {\tt grabiner@wcnet.org}
\end{tabular}
\end{document}